\documentclass[ijaa,sglanonrev]{informs4}

\makeatletter
\def\theARTICLETOP{}
\makeatother

\usepackage{eqndefns-left} 
\RequirePackage{tgtermes}
\RequirePackage{newtxtext}
\RequirePackage{newtxmath}
\RequirePackage{bm}
\RequirePackage{endnotes}

\OneAndAHalfSpacedXII 

\usepackage{algorithm}
\usepackage{algpseudocode}
\usepackage{tikz}
\usepackage{booktabs}
\usepackage[left]{lineno} 
\usepackage{placeins}
\usepackage{xcolor}

\usepackage{xcolor, tabularx}

\usepackage{makecell}

\definecolor{darkgreen}{RGB}{20,120,70} 

\usepackage{natbib}
 \bibpunct[, ]{(}{)}{,}{a}{}{,}%
 \def\bibfont{\small}%
\EquationsNumberedThrough    

\TheoremsNumberedThrough     
\ECRepeatTheorems  %

\MANUSCRIPTNO{IJAA-0001-2024.00}

\begin{document}


\RUNAUTHOR{Kilb, Bienstock and Newman}

\RUNTITLE{Humanitarian Logistics under Uncertainty}

\TITLE{Probabilistic Modeling versus Robust Optimization: A tutorial based on a humanitarian logistics use case}

\ARTICLEAUTHORS{%
\AUTHOR{Justin Kilb}
\AFF{Operations Research with Engineering Graduate Program, Mechanical Engineering Department, Colorado School of Mines, Golden CO 80401 \EMAIL{justinkilb@mines.edu}}

\AUTHOR{Daniel Bienstock}
\AFF{Industrial Engineering and Operations Research Department, Columbia University, New York, NY 10027
\EMAIL{dano@columbia.edu}}

\AUTHOR{Alexandra M. Newman}
\AFF{Operations Research with Engineering Graduate Program, Mechanical Engineering Department, Colorado School of Mines, Golden CO 80401
\EMAIL{anewman@mines.edu}}
} 

\ABSTRACT{%
This tutorial contrasts probabilistic modeling and robust optimization to determine  decisions in humanitarian logistics, specifically supply chains subject to adversarial (natural and human) disruptions. Natural disruptions induce dispatch  of long-haul relief supply movement  as storm forecasts evolve. A two-step workflow: (i) computes an initial pre-staging plan from the most likely forecast, and (ii) evaluates that fixed plan across plausible deviations in the eventual landfall location. In this way, dispatch decisions balance lead time and improved forecast information. For last-mile distribution, we propose deliveries when transportation networks must be protected against the worst case. We apply an iterative robust routing method that detects high-concentration links and increases their effective cost to promote route diversification. A case study based on Typhoon Noru (2022)  shows how the combined approach identifies an optimal dispatch time and then protects last-mile delivery from difficult-to-predict network disruptions that could jeopardize the entire supply-chain operation.
}%

\FUNDING{This research was supported by Office of Naval Research Grant Numbers N000142412442, N0001424B001, and by a generous donation from Defcon AI.}




\maketitle

\section{Introduction}
Practitioners use operations research to support humanitarian decision making in preparedness and response, with an emphasis on translating analytic models into measurable societal impact \citep{Altay2006, ErgunKeskinocakSwann2011, Peters2022}. Humanitarian logistics requires rapid delivery of essential aid in a non-deterministic supply chain, which differs fundamentally from that typical of a commercial system due to the degree of uncertainty  and significant infrastructure disruptions  \citep{VanWassenhove2006, Hu2023}. Both probabilistic and robust optimization frameworks manage these complexities.

Probabilistic optimization represents uncertainty invoking distributions or scenario sets and optimizes an expected performance measure or another probability-based criterion. When decision makers cannot justify reliable probabilities or prefer protection against worst-case realizations, robust optimization is the method of choice.  Generally, this paradigm specifies an uncertainty set,  selecting solutions that satisfy the model constraints for all realizations in that set. Performance is optimized under the most unfavorable realizations in the set. Both approaches support resilient humanitarian planning, but organizations often deploy them under different degrees of uncertainty and risk tolerances.

This tutorial synthesizes these two perspectives to address critical humanitarian dispatch decisions for settings in which an incoming storm event (e.g., a hurricane or typhoon) is forecastable but uncertain, and predictions evolve as the event approaches. We first examine the problem of \textbf{dispatch timing} through a probabilistic lens, using a two-phase workflow denoted ($\mathcal{P}_{\!1}$) and ($\mathcal{P}_{\!2}$), respectively. In this context, planners must decide when to move relief supplies from distant international hubs to regional staging areas. This decision involves a fundamental trade-off: dispatching early can accelerate delivery but increases the risk of mispositioning or losing supplies if the storm’s realized path deviates, whereas dispatching later benefits from improved forecast accuracy at the cost of reduced transportation lead time.

We then transition to the \textbf{last-mile distribution} problem using a robust optimization approach. Once supplies reach regional hubs, the primary uncertainty shifts from the event's location to network accessibility. Here, robust methods protect against difficult-to-predict localized infrastructure failures, helping to ensure that the pre-staged supplies reach affected populations even when transportation corridors are compromised. We denote the baseline \emph{risk-unaware} plan by ($\mathcal{R}^{-}$) and the \emph{robust} alternative by ($\mathcal{R}^{+}$).

This tutorial elucidates how to:
\begin{enumerate}
\item Select a probabilistic or robust optimization approach based on the decision context, the quality of available uncertainty information, and the decision-maker’s risk tolerance.
\item Construct a forecast-driven dispatch timing workflow that maps evolving forecasts to time-indexed expectations of closing time and unmet demand. 
\item Implement an iterative robust last-mile routing method that identifies vulnerable high-flow links and diversifies routes to protect delivery performance under localized network disruptions.

\end{enumerate}

We first use probabilistic optimization to show how planners choose when to dispatch supplies under evolving forecast uncertainty, placing this decision in the broader context of stochastic pre-positioning models.

\section{Stochastic optimization in humanitarian dispatch}


Foundational humanitarian logistics applications use two-stage stochastic programs for strategic prepositioning. In the first stage, planners select facility locations and inventory levels before the disaster outcome is known. In the second stage, after demand and network conditions realize, they allocate supplies and route deliveries through recourse decisions. \cite{Rawls2010} develop a stochastic mixed-integer program  to locate emergency supplies in anticipation of hurricanes, accounting for uncertain demand and network availability. \cite{Balcik2008} integrate facility location and inventory pre-positioning decisions in a maximal-covering model to maximize expected covered demand, and \cite{Duran2011} demonstrate the benefits of centralized planning for international relief items. These models reduce expected response times but fix the time of supply dispatch from distant hub to the affected region in advance. They optimally stage supplies for a given forecast snapshot, rather than optimizing when to dispatch as forecasts evolve.

To capture the dynamic nature of information arrival, researchers extend the two-stage framework to three or more stages. These formulations allow sequential decisions to be revised across multiple, discrete time periods as situational awareness improves, explicitly modeling the value of information throughout the forecast horizon. \citet{Salmeron2011} introduce the Global Fleet Station Mission Planner, an optimization-based mission planning and scheduling tool that supports rapid rescheduling in response to exigencies and changing operational requirements. \cite{TAVANA201821} develop a multi-period, mixed-integer optimization framework that jointly plans pre-disaster warehouse siting and perishable inventory decisions and then, conditional on when an earthquake occurs, optimizes post-disaster two-echelon vehicle routes to distribute relief to affected areas. \citet{Hu2023} formulate dynamic relief pre-positioning as a multi-stage stochastic program in which procurement (and inventory via procurement and returns) is adapted over time as uncertainty unfolds along a scenario tree. Recent applied work adopts scenario-based, rolling-horizon optimization to study end-to-end humanitarian supply chains, jointly examining advance demand information, prepositioning strategies, and last-mile distribution performance using real operational data \citep{Jing2025}.

Although formulations with additional stages capture the value of waiting for more accurate information, they typically encode the timing of initial supply dispatch through a stochastic stage structure and scenario tree, rather than treating it as an explicit decision variable. A notable exception is \cite{Rezapour2021}, who directly study dispatch timing by solving a stochastic prepositioning model across a range of trigger times and then by identifying the point at which logistics cost and response time are balanced. Their results highlight the fundamental trade-off between improved forecast accuracy and the operational costs of delayed action. However, their approach treats the time at which supplies are dispatched as an exogenous parameter that is later varied using post-optimality analysis, rather than as an endogenous decision. Therefore, the underlying scenario structure does not explicitly represent how forecast uncertainty evolves through intermediate threat (in this case, storm) positions.

\begin{figure}[!htb]
\FIGURE
{\includegraphics[width=\linewidth]{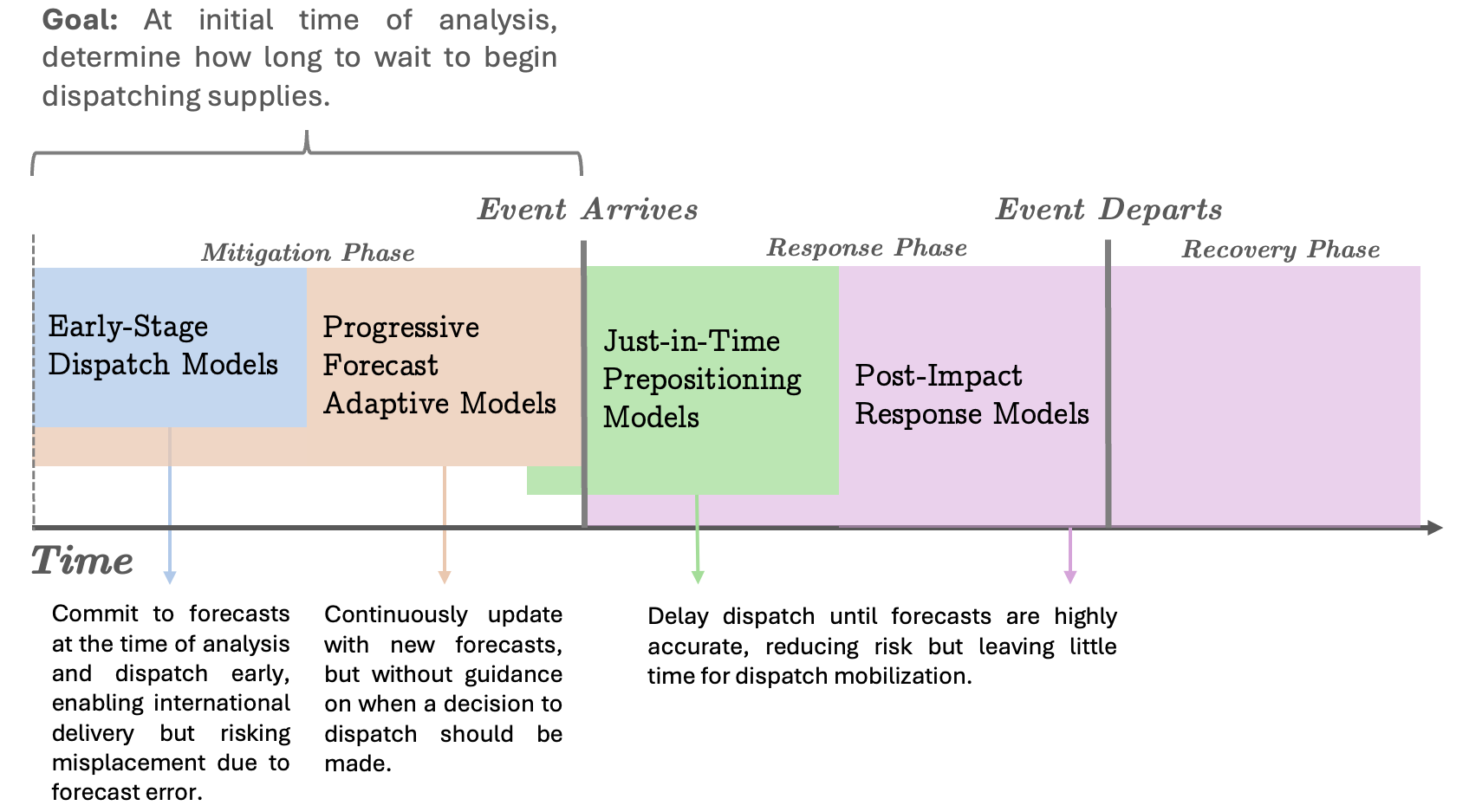}}
{Timeline summary of dispatch and prepositioning model classes and the decision point addressed in this tutorial.
\label{fig:lit_review_timeline}}
{}
\end{figure}

Existing stochastic models that support dispatch and prepositioning also differ in when they place the decision relative to the event timeline. Figure~\ref{fig:lit_review_timeline} summarizes four common categories: (i) early-stage dispatch models that commit at the initial time of analysis, (ii) progressive forecast adaptive models that update plans as new forecasts arrive but do not explicitly identify when to commit, (iii) just-in-time prepositioning models that delay dispatch until forecasts tighten near event arrival, and (iv) post-impact models that allocate and route supplies after the event occurs. This tutorial focuses on a progressive forecast adaptive model (ii) that occurs within the mitigation phase. 

This tutorial builds on these concepts in the probabilistic modeling domain while introducing three methodological distinctions. First, dispatch time is modeled as an explicit decision variable optimized within the framework, rather than selected through post-hoc comparison. Second, the method represents forecast evolution through a structured decomposition of intermediate event positions and conditional outcome probabilities, which makes the narrowing of the forecast uncertainty region directly observable. Third, the two-phase architecture separates spatial allocation from temporal optimization, providing a modular and interpretable workflow that allows practitioners to examine each layer of the decision process independently. This transparency makes the approach well suited to the tutorial format and accessible to readers seeking an introduction to forecast-driven dispatch optimization without the complexity of monolithic stochastic programs.

\section{Methodology | Probabilistic optimization for dispatch timing}\label{sec:Method}

This section introduces a two-phase modeling framework designed for operational utility and transparency. The approach separates the decision process into two sequential phases to reflect the planning structures of humanitarian organizations. Model ($\mathcal{P}_{\!1}$) addresses the spatial challenge by selecting pre-stage locations and allocations that minimize delivery time while accounting for the risk of site destruction; ($\mathcal{P}_{\!2}$) solves the temporal problem by utilizing ($\mathcal{P}_{\!1}$) outputs to determine the optimal dispatch time. This process balances the advantages of improved forecast accuracy against the operational cost of lost lead time.

\subsection{Problem setup}
\label{sec:problem_setup}

To formalize the dispatch timing problem, we define the logistics system elements and the assumptions linking the system to forecast uncertainty. 

\subsubsection{System description}

\begin{itemize}
    \item \textbf{Dispatch location:} The distant source from which supplies are initially shipped. Long-haul transport durations create the cost that necessitates early action.
    \item \textbf{Pre-stage locations:} Candidate sites near the impact zone used for temporary storage to reduce post-event delivery time.
    \item \textbf{Closing time:} The total duration required to complete all deliveries from pre-stage sites to demand nodes; this represents the primary measure of response effectiveness.
    \item \textbf{Dispatch lead time:} The duration required to move all supplies from the source to pre-stage locations; this  directly influences the trade-off against improving forecast accuracy.
    \item \textbf{Demand locations:} Areas within the projected event trajectory inferred from forecasts which may shift as new information arrives.
\end{itemize}

\subsubsection{Forecast modeling assumptions}

The framework uses externally provided probabilistic forecasts as model inputs, which ensures compatibility with diverse forecasting standards. The framework assumes that forecast error grows as the   horizon lengthens, so forecasts issued farther from event arrival contain greater spatial uncertainty than those issued closer to arrival. This pattern is consistent with dynamical forecasting: small errors in initial conditions and model representation can amplify over time. This assumption matters because dispatch lead time determines how early planners must commit. Short-haul deployment permits waiting for late forecasts to unfold before demand materializes. Long-haul international deployment requires commitment while forecast uncertainty remains large, so dispatch timing governs the trade-off between earlier action and improved information. For example, if organizations can dispatch and deploy supplies within hours, they can often wait for later forecasts and still meet response timelines. If deployment requires days, planners might  dispatch even while forecasts   contain substantial uncertainty.

\subsection{Phase 1 ($\mathcal{P}_{\!1}$): Pre-stage location optimization}

Model ($\mathcal{P}_{\!1}$) addresses the spatial allocation of supplies to minimize delivery time while mitigating the risk of site-level disruptions and uses the most likely event trajectory from the initial forecast as the baseline landfall location for planning. The following sets define the decision space for the optimization:

\begin{table}[h!]
	\begin{tabular}{l l}
		{\bf Symbol} & {\bf Definition} \\
		$\mathcal{I}$ & candidate pre-stage locations \\
		$\mathcal{J}$ & anticipated demand locations \\
		$\mathcal{K}$ & supply types (e.g., water, shelter) \\
		$\mathcal{M}$ & transportation modes (e.g., truck1, truck 2) \\
	\end{tabular}
\end{table}

The ($\mathcal{P}_{\!1}$) model determines the quantity of each supply type $k$ to allocate to each candidate pre-stage location $i$ for delivery to demand location $j$ using vehicle mode $m$. To represent the potential loss of supplies, we assign each pre-stage location $i$ a probability $p_i$ of becoming inaccessible or being destroyed based on exposure to the event and also based on the site resilience. For example, the pre-stage location might be in the event trajectory, but underground in a bunker, in which case the user might assign a lower $p_i$. This probability applies uniformly to all commodities stored at a specific location.

The parameter $\alpha_k$ weights the disruption-risk penalty relative to closing time in the objective, allowing the model to trade delivery efficiency against event exposure separately for each supply type $k$. Because the disruption term is evaluated at every candidate pre-staging site $i$, $\alpha_k$ effectively controls how strongly the model prefers moving supply type $k$ to locations with lower disruption probability. Larger values of $\alpha_k$ represent lower tolerance for exposure and therefore push the pre-staging allocation for supply type $k$ away from the forecasted trajectory, while smaller values permit greater reliance on sites closer to the projected impact region.

\begin{figure}[!htb]
\FIGURE
{\includegraphics[width=\linewidth]{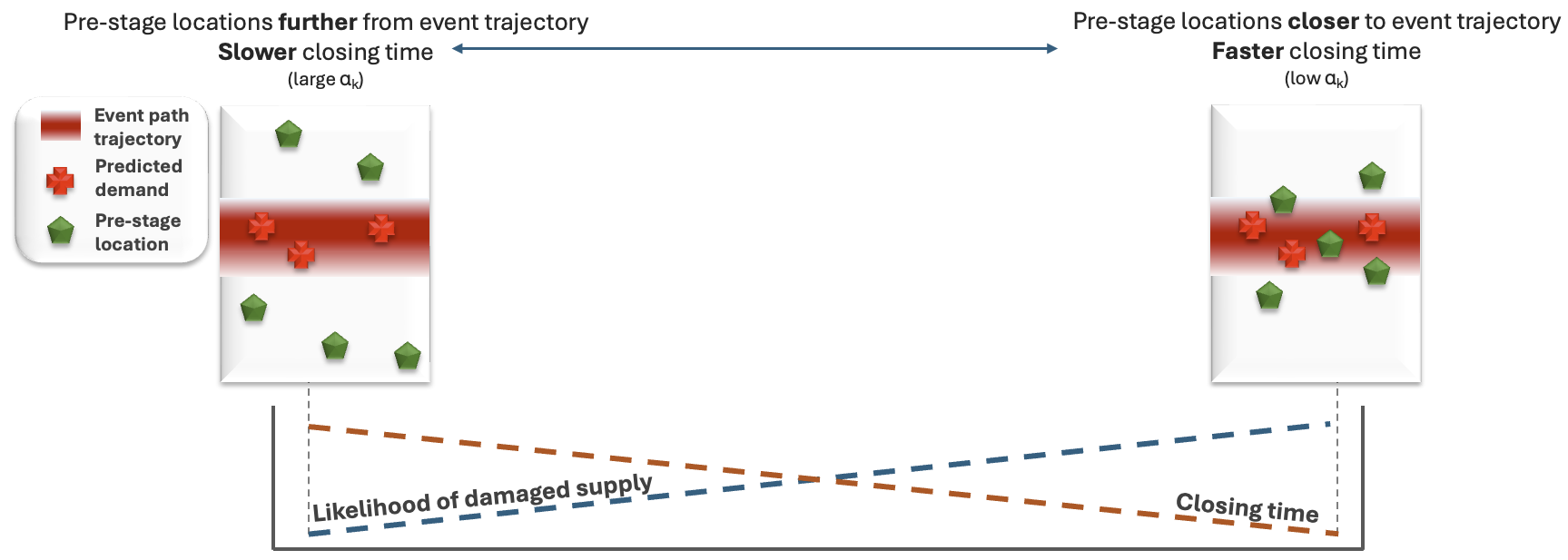}}
{Effect of increasing $\alpha_k$ with ($\mathcal{P}_{\!1}$). Higher values increase the penalty associated with storing supplies at locations inside the current event trajectory, thereby shifting the initial pre-stage allocation outward.
\label{fig:alpha_term}}
{}
\end{figure}

Supply-specific tuning is important because different supplies face different operational constraints. For example, some items require a fast closing time but are relatively low-cost and small in volume and weight. In such cases, planners may accept staging a portion of supply type $k$ near the projected trajectory to pursue rapid delivery, while simultaneously staging backup inventory at safer sites farther from the event so that service can continue if exposed stock is lost. Importantly, the additive penalty does not impose a hard constraint that forbids using exposed sites; instead, it allows them to be used sparingly when they materially improve delivery time, which can be especially valuable in sparse networks, where there are only a few feasible pre-staging sites and excluding nearby options can force inventory to be positioned at  locations   much farther away. Figure~\ref{fig:alpha_term} illustrates how larger values of $\alpha_k$ shift pre-staging allocations away from the forecasted trajectory.

The ($\mathcal{P}_{\!1}$) model addresses the spatial allocation of supplies to minimize closing time while mitigating the risk of site-level disruptions and uses the most likely event trajectory from the initial forecast as the baseline landfall location for planning. 

\subsection{From a single forecast to scenario-based evaluation}
\label{sec:phase1-to-scenarios}

The optimal pre-staging plan generated in ($\mathcal{P}_{\!1}$) assumes the forecasted trajectory is certain. However, plans optimized for the most likely track may perform poorly if event landfall locations shift. This motivates a structured extension to evaluate plan performance across alternative event trajectories. To assist decision-makers in quantifying the operational consequences of forecast error, we introduce a scenario-based evaluation method. Unlike traditional stochastic approaches that embed uncertainty within a single optimization problem, this framework separates these layers to ensure forecast deviation effects remain observable and comparable. By externalizing performance degradation, the framework provides a clear understanding of how exposed a plan  becomes when reality differs from initial expectations.

\subsubsection{Scenario and outcome sensitivity analysis}

Forecast uncertainty is represented through a set of discretized landfall outcomes $\mathcal{O} = \{1,2,\dots,|\mathcal{O}|\}$. Outcome 1 denotes the most likely event trajectory used in ($\mathcal{P}_{\!1}$), and the remaining outcomes represent possible forecasted deviations on either side of the most likely trajectory. The most likely outcome, outcome 1, is determined by the forecast chosen by the user at the time of analysis. The required number of outcomes depends on forecast spread and on how sensitively the logistics model responds to event landfall shifts. The forecast cone (or ``cone of uncertainty'') summarizes plausible future storm-center trajectories and typically widens farther into the future. Wide forecast cones and dense sets of candidate pre-stage and demand locations often warrant more outcomes because small changes in landfall can alter the exposed sites and the induced demand pattern, producing materially different closing times and unmet-demand levels. By contrast, when the network and demand geometry are coarse relative to the forecast cone, many nearby landfall event trajectory realizations yield operationally similar consequences, so fewer outcomes can suffice. This tutorial uses five outcomes to balance interpretability with geographic representation in the ($\mathcal{P}_{\!1}$) post-processing step.

\subsubsection{($\mathcal{P}_{\!1}$) post-processing outcome sensitivity analysis}
The post-processing stage evaluates the fixed ($\mathcal{P}_{\!1}$) pre-staging solution across all deviated forecast outcomes $o \in \mathcal{O}$. The solution provided by ($\mathcal{P}_{\!1}$) fixes the staged quantities $X_{ijkm}$ for each supply type $k$ at location $i$ assigned to demand node $j$ using mode $m$. We treat these $X_{ijkm}$ values as committed decisions in order to interpret the consequences of adopting this pre-staging strategy under each deviated outcome. Operationally, the model behaves as if the pre-staged supplies have already arrived at their assigned locations, so subsequent solutions reflect performance given this realized positioning strategy.

For each outcome $o$, we re-solve the distribution problem with the pre-staged inventory fixed, i.e., with the total quantity of each supply type $k$ stored at each pre-stage site $i$ held constant. The model can reassign those fixed stocks across new demand nodes $j$ and transportation modes $m$ to reflect the realized event trajectory and the resulting demand pattern, but it cannot relocate inventory across pre-stage sites. If a pre-stage location falls within the impact zone of outcome $o$, its inventory is considered inaccessible. 

This analysis produces two parameters for each supply type $k$ under outcome $o$: closing time $c_{ko}$ and percent unmet demand $\bar{d}_{ko}$. Given the probability $\check{p}_{o}$ of outcome $o$ occurring, expected scenario performance is defined as: 

\[
\begin{aligned}
e_{k} &= \sum_{o \in \mathcal{O}} \check{p}_{o}\, c_{ko} \quad \forall k \in \mathcal{K}\\
\bar{d}_{k} &= \sum_{o \in \mathcal{O}} \check{p}_{o}\, \bar{d}_{ko} \quad \forall k \in \mathcal{K}
\end{aligned}
\]

\begin{figure}[!htb]
\FIGURE
{\includegraphics[width=0.70\linewidth]{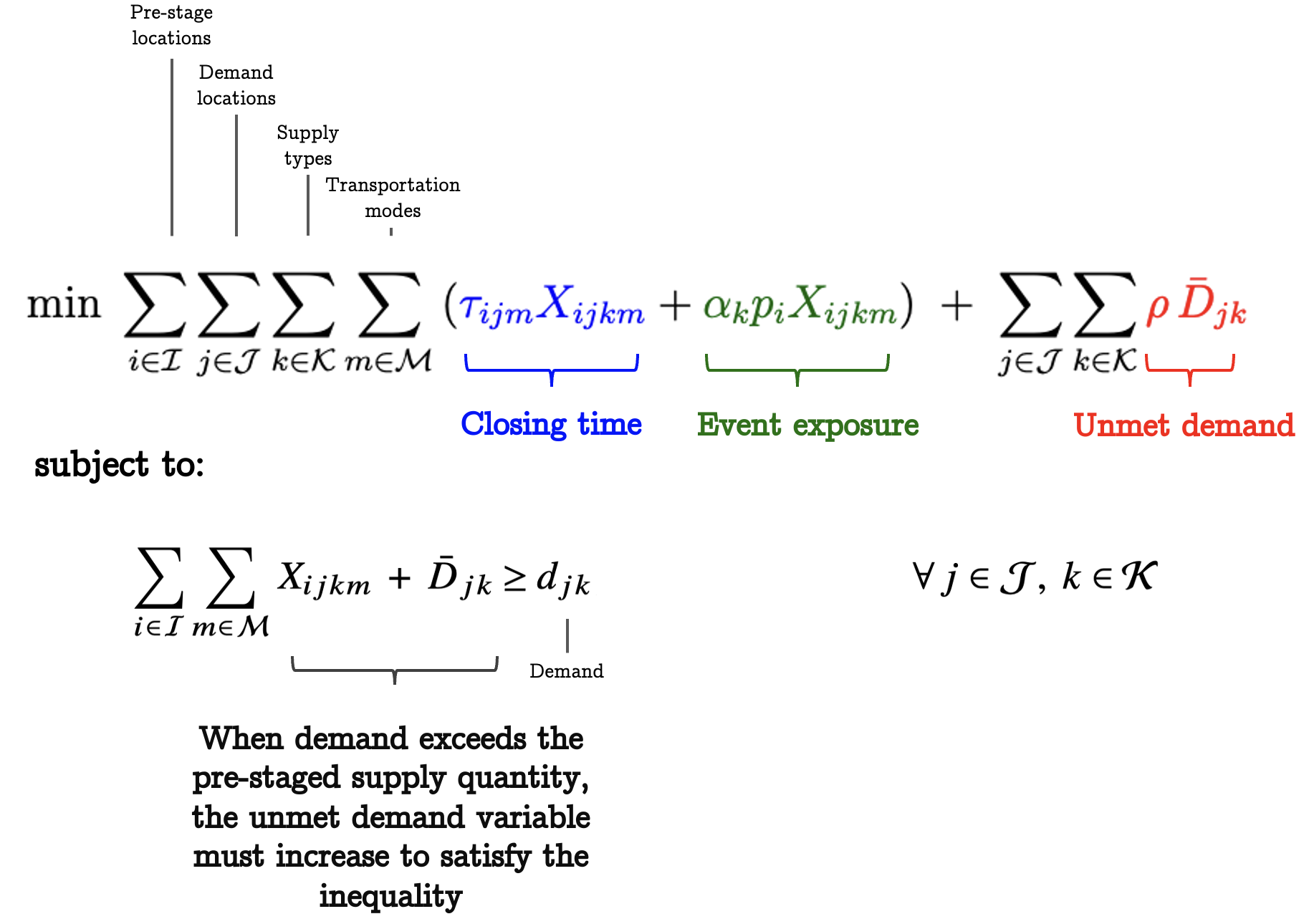}}
{Objective function and representative constraint with unmet demand. The figure introduces $\bar{D}_{jk}$ to capture shortfalls when inventory becomes inaccessible or demand increases under alternative-event landfall outcomes.
\label{fig:phase1_with_demand}}
{}
\end{figure}

To track shortages created by inventory loss or demand surges under deviated outcomes, we introduce the unmet demand variable $\bar{D}_{jk}$. Figure~\ref{fig:phase1_with_demand} provides the first view of the formulation by showing the objective function and a representative demand-balance constraint where capital letters in Roman font indicate a decision variable. The full formulation can be found in the appendix. 

\begin{figure}[H]
\FIGURE
{\includegraphics[width=\linewidth]{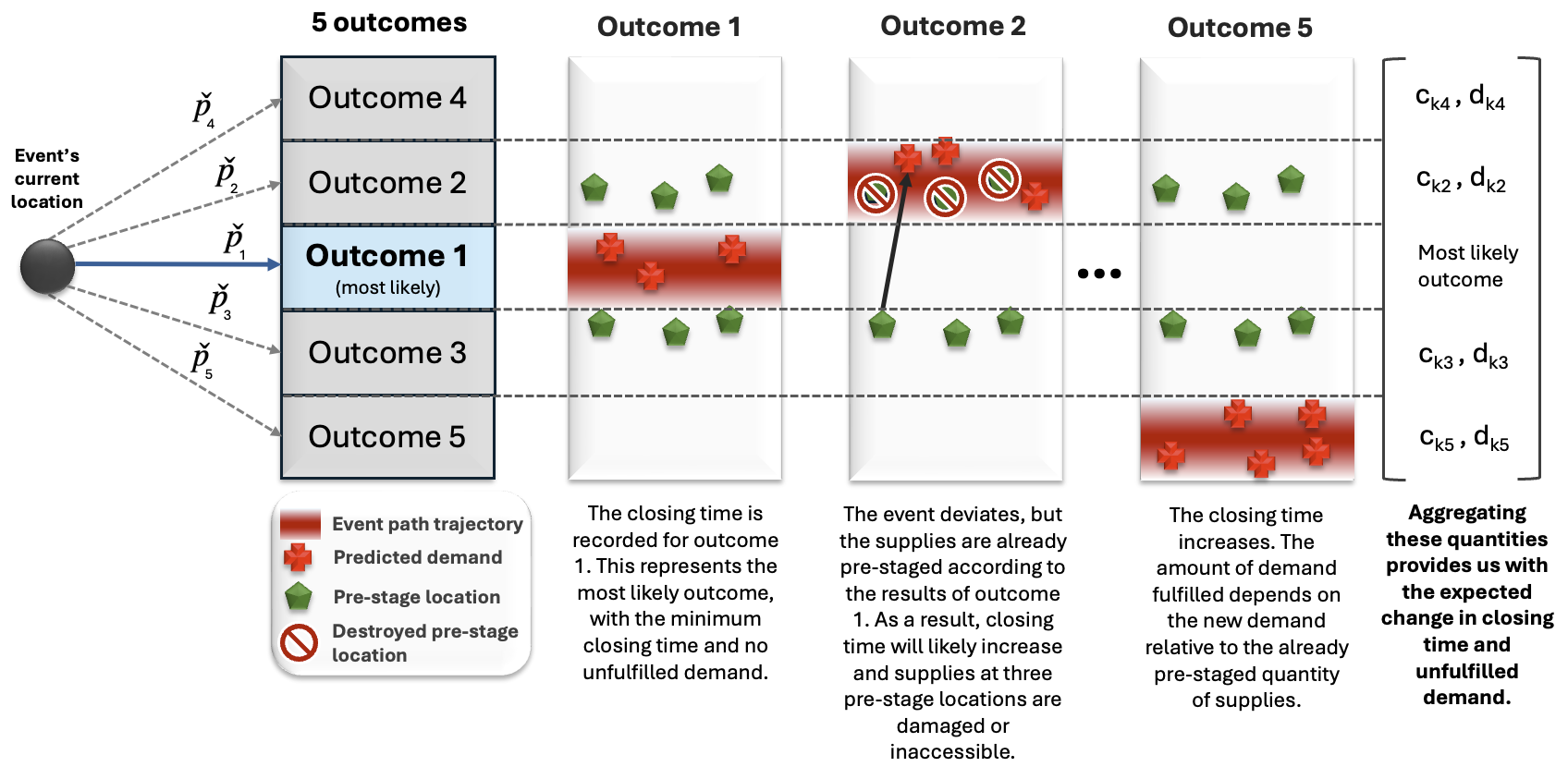}}
{Post-processing evaluation across event landfall outcomes. The figure shows how forecast deviations affect demand fulfillment and closing time when the model pre-stages supplies using outcome 1.
\label{fig:scenario_robustness}}
{}
\end{figure}

Figure~\ref{fig:scenario_robustness} summarizes the post-processing step by evaluating the fixed pre-staging plan across outcomes, with outcome 1 as the most likely forecast.

\subsubsection{Discretizing forecasted intermediate event locations}
The metrics described previously are extended to incorporate an event forecast evolving through intermediate forecast locations. This extension evaluates the consequences of dispatch timing as the storm progresses and the forecast cone narrows. We introduce the set of intermediate forecast positions $\mathcal{L} = \{1,2,\dots,|\mathcal{L}|\}$ and index results by scenario $s \in \mathcal{S}$, which corresponds to a specific forecast horizon. The intermediate event locations can be selected directly from operational forecasts issued during an event. Forecasters typically assign probabilities to striking locations along its trajectory toward  eventual landfall, along with a forecast cone that widens farther into the future. As the event advances toward landfall, forecast uncertainty typically decreases and the implied uncertainty over the landfall location narrows; this information refinement is the mechanism we leverage when evaluating dispatch timing. This tutorial does not address how to construct or calibrate these probabilities. Instead, we treat the forecast track and associated uncertainty as inputs, so the same procedure applies to any approaching event in which trajectory accuracy improves and the uncertainty cone contracts as the event nears the region of impact. 

Performance metrics are modified to include this event location and scenario index. The parameter $c_{kos}$ represents the closing time, and $\check{d}_{kos}$ represents the percent unmet demand for all $k$, $o$, and $s$. At each event location $\ell$, we track outcome uncertainty $\check{p}_{o\mid \ell,s}$, the probability of outcome $o$ conditional on storm position $\ell$ and scenario $s$, and event trajectory uncertainty $\hat{p}_{\ell\mid t}$, the probability that the storm occupies location $\ell$ conditional on time $t \in \mathcal{T} = \{1,2,\dots,|\mathcal{T}|\}$, where the elements of $\mathcal{T}$ index discrete time steps preceding event landfall. Throughout, the expected metrics are taken with respect to the discrete probability weights defined above. Equivalently, these are conditional expectations under the available forecast information, but we suppress additional conditioning notation for ease of exposition. For example, $\hat{p}_{\ell\mid t}$ is shorthand for the more explicit conditional probability $\mathbb{P}(L_t=\ell \mid \mathcal{F}_t)$, where $L_t$ denotes the (discretized) storm location at time $t$ and $\mathcal{F}_t$ denotes the forecast at time $t$. In keeping with the tutorial format, we employ index-based notation and suppress additional conditioning to maintain clarity throughout the formulation.

\begin{figure}[H]
\FIGURE
{\includegraphics[width=\linewidth]{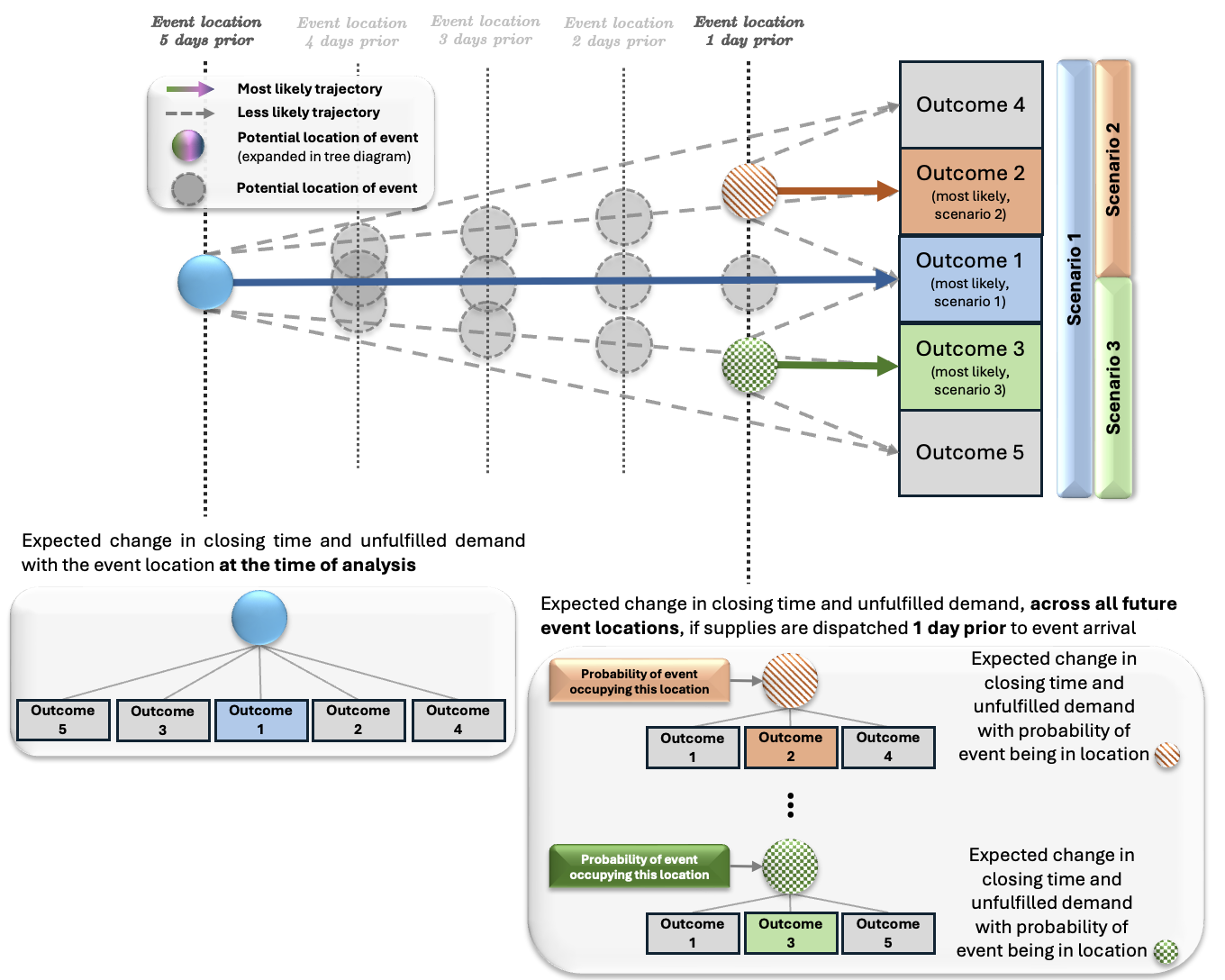}}
{Scenario and outcome structure used to evaluate the ($\mathcal{P}_{\!1}$) pre-staging plan under evolving forecast uncertainty.
\label{fig:phase2}}
{}
\end{figure}

We compute the expected performance for each supply type $k$ at each location $\ell$:

\[
\begin{aligned}
e_{k\ell} &= \sum_{s \in \mathcal{S}} \sum_{o \in \mathcal{O}} \check{p}_{o\mid \ell,s}\, c_{k o s} \quad \forall k \in \mathcal{K},\ \forall \ell \in \mathcal{L}\\
\bar{d}_{k\ell} &= \sum_{s \in \mathcal{S}} \sum_{o \in \mathcal{O}} \check{p}_{o\mid \ell,s}\, \check{d}_{k o s} \quad \forall k \in \mathcal{K},\ \forall \ell \in \mathcal{L}
\end{aligned}
\]
Intermediate locations reshape the probability distribution over outcomes without creating additional optimization stages. Using these metrics, we aggregate across event locations to compute dispatch-time expectations:
\[
\begin{aligned}
\hat{e}_{kt} &= \sum_{\ell \in \mathcal{L}} \hat{p}_{\ell\mid t}\, e_{k\ell} \quad \forall k \in \mathcal{K},\ \forall t \in \mathcal{T}\\
\hat{d}_{kt} &= \sum_{\ell \in \mathcal{L}} \hat{p}_{\ell\mid t}\, \bar{d}_{k\ell} \quad \forall k \in \mathcal{K},\ \forall t \in \mathcal{T}
\end{aligned}
\]
These expectations quantify the predicted change in system performance if supplies are dispatched at time $t$. Figure~\ref{fig:phase2} illustrates this shift in perspective as the model transitions from outcome-level expectations to the time-dependent expectations required for ($\mathcal{P}_{\!2}$). 

\subsection{Dispatch timing optimization: ($\mathcal{P}_{\!2}$)}

The ($\mathcal{P}_{\!2}$) model utilizes results from ($\mathcal{P}_{\!1}$) and the post-processing stage to evaluate the trade-off between early dispatch under uncertain forecasts and late dispatch under improved accuracy with reduced lead time. For each supply type $k$, the model selects a single dispatch time from the candidate set $\mathcal{T}$ using the binary decision variable $Z_{kt}$, which equals $1$ if the model dispatches supply type $k$ at time $t$ and equals $0$ otherwise. The ($\mathcal{P}_{\!2}$) model then evaluates the expected increase in closing time resulting from delayed mobilization and the expected percentage of unmet demand.

The first objective term penalizes delayed delivery by comparing the fixed lead time from the source to pre-stage locations with the expected closing time $\hat{e}_{kt}$ derived from scenario-outcome analysis. The second term penalizes expected unmet demand $\hat{d}_{kt}$. User-defined weights $a_k$ and $b_k$ represent the relative importance of closing time and unmet demand for each supply type $k$. These weights also serve as scaling factors to account for differing units of measurement, which prevents either metric from dominating the objective purely due to magnitude.

\begin{figure}[H]
\FIGURE
{\includegraphics[width=.7\linewidth]{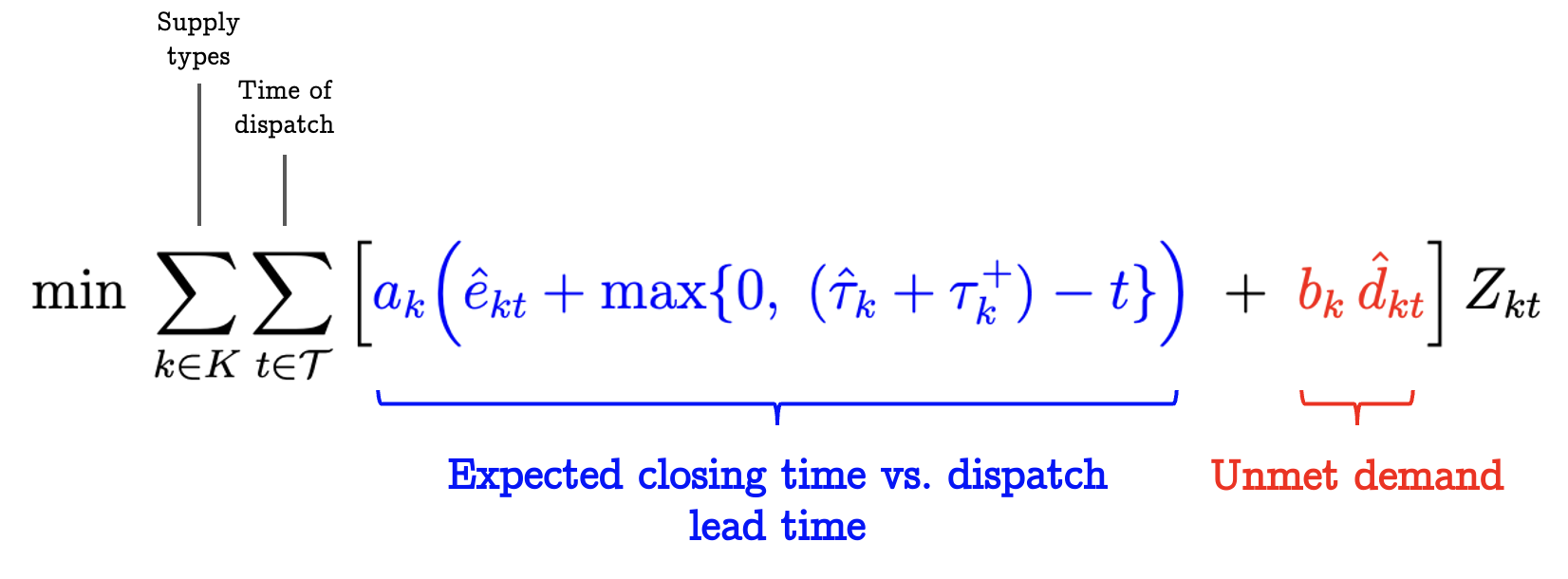}}
{The first objective term penalizes expected closing time and any lateness beyond the available lead-time buffer, while the second term penalizes the expected fraction of unmet demand. The binary variable $Z_{kt}$ selects exactly one dispatch time for each supply type.
\label{fig:phase_2_obj}}
{}
\end{figure}

A delay factor $\tau_k^+$ provides a temporal buffer by postponing dispatch according to practitioner risk tolerance. In effect, $\tau_k^+$ can override the objective tradeoff for specific supplies: a user may prefer to dispatch a given supply type later, accepting delayed closing time, provided the plan guarantees successful delivery under the realized forecast. This is especially relevant for commodities that are not highly time sensitive but are expensive (or otherwise costly to lose), for which decision makers may be unwilling to stage early under uncertainty. Figure~\ref{fig:phase_2_obj} summarizes the resulting ($\mathcal{P}_{\!2}$) objective function. The term $\max\{0,\,(\hat{\tau}_k + \tau_k^+) - t\}$ imposes a lateness penalty when the chosen dispatch time provides insufficient lead time. This operator is linearized in the Appendix via an auxiliary variable and inequality constraints.

\subsection{Summary of dispatch timing methodology}

The ($\mathcal{P}_{\!1}$) model utilizes a single-event forecast, candidate locations, and transport modes to determine the optimal pre-staging plan $X_{ijkm}$. The model balances proximity to demand against the risk of site failure using a location-specific risk term and a supply-type tuning parameter. The resulting output defines the spatial allocation and movement strategies alongside performance measures for the most likely forecast trajectory.

The ($\mathcal{P}_{\!1}$) post-processing and ($\mathcal{P}_{\!2}$) model evaluates the temporal trade-off between waiting for forecast improvements and the resulting loss of dispatch lead time. By assessing the fixed pre-staging plan against shifted event trajectories, the method evaluates expected changes in closing time $\hat{e}_{kt}$ and unmet demand $\hat{d}_{kt}$. These metrics, combined with dispatch lead times, enable the selection of an optimal dispatch time for each supply type.

To clarify the methodological position for readers familiar with stochastic optimization, ($\mathcal{P}_{\!1}$) is a facility-location model with a risk-weighted objective optimized under a single forecast. Uncertainty is managed in a separate scenario sensitivity analysis in which $\hat{e}_{kt}$ and $\hat{d}_{kt}$ are computed via exact summation across an explicitly enumerated scenario tree. These expectations serve as fixed parameters in ($\mathcal{P}_{\!2}$), rather than as random variables.

Consequently, this framework is a two-phase deterministic-stochastic hybrid, rather than a classical two-stage stochastic program. In standard notation, a two-stage model includes recourse variables indexed by the realized outcome $\omega \in \Omega$ to adapt decisions after uncertainty resolves. This framework includes no outcome-indexed recourse decisions because it solves both phases at the initial planning time. The ($\mathcal{P}_{\!2}$) model instead selects a dispatch time that determines when the organization commits to mobilization. 

The distinguishing feature of this methodology is the explicit modeling of forecast-cone evolution. By structuring uncertainty through outcome sets $o$, intermediate locations $\ell$, and conditional probabilities $\check{p}_{o\mid \ell,s}$ and $\hat{p}_{\ell\mid t}$, the approach constructs time-dependent performance curves. This transparent workflow remains interpretable while providing a modular foundation for advanced extensions such as the incorporation of continuous distributions. 

From a computational perspective, the framework remains tractable because it does not embed the scenario tree inside a single large optimization model. The mathematical representations of ($\mathcal{P}_{\!1}$) and ($\mathcal{P}_{\!2}$) remain parsimonious. The framework computes scenario-based performance measures in a separate post-processing step by re-solving the distribution model across outcomes with the ($\mathcal{P}_{\!1}$) plan fixed.

\subsection{Case study | Probabilistic optimization for dispatch timing}\label{sec:case_study}

We apply this two-phase framework  to a historical event characterized by high forecast uncertainty and significant dispatch lead times. The following application demonstrates the integration of pre-staging, scenario-based evaluation, and timing optimization using actual forecast evolution and geographic exposure data.

\subsection{Typhoon Noru as an example}

Typhoon Noru, locally designated as Super Typhoon Karding, formed on September 21, 2022. The storm underwent rapid intensification before making its first landfall over the Polillo Islands at 07:30 Coordinated Universal Time on September 25. The persistence of forecast uncertainty just days before landfall, coupled with the extended lead times required for international aid, provides a suitable context for evaluating the trade-off between forecast accuracy and early mobilization. The case study utilizes three primary data components comprising probabilistic event forecasts, estimated humanitarian demand, and a logistics network with associated travel times.

\subsubsection{Typhoon forecast data}

We obtain event forecast data from the European Centre for Medium-Range Weather Forecasts \citep{ecmwf_datasets}. Importantly, each forecast we use is ex ante: it reflects the information available at the time of analysis and represents the provider's prediction of the event's future locations together with its forecast cone. For this case study, we use a strike-probability map reported in six-hour increments, which provides the likelihood that the storm center will pass within 120~km of each grid cell. We interpret these values as relative likelihoods and normalize them to obtain weights that sum to one over the discretized set of outcomes in each scenario.

\begin{figure}[H]
\FIGURE
{\includegraphics[width=0.60\linewidth]{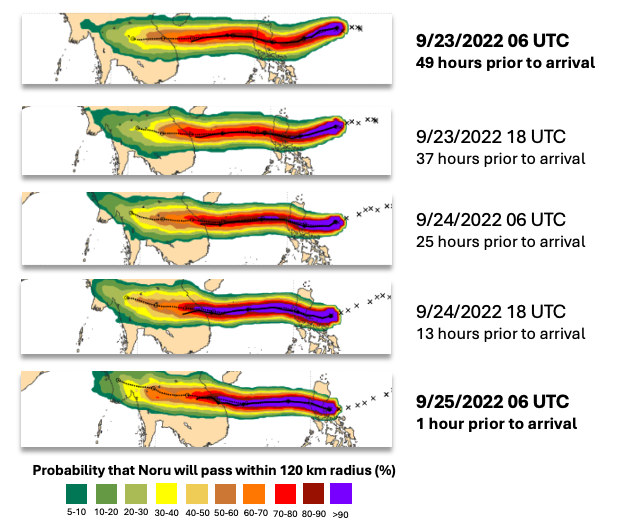}}
{Evolution of strike probability forecasts for Typhoon Noru showing the reduction in geographic uncertainty. \label{fig:noru_forecasts}}
{Source: Adapted from European Centre for Medium-Range Weather Forecasts data.}
\end{figure}

We do not claim that this is the only way to obtain probabilities of event location or to quantify forecast uncertainty; forecasting methods and uncertainty quantification are active and evolving research areas. Our framework requires only that, at each forecast issuance time, the user can provide: (i) a probabilistic description of future event locations; and, (ii) an uncertainty representation over deviations in the eventual impact region (e.g., a forecast cone or an equivalent probabilistic summary). Figure~\ref{fig:noru_forecasts} shows the strike probability maps used in this case study.

\begin{figure}[!htb]
\FIGURE
{\includegraphics[width=\linewidth]{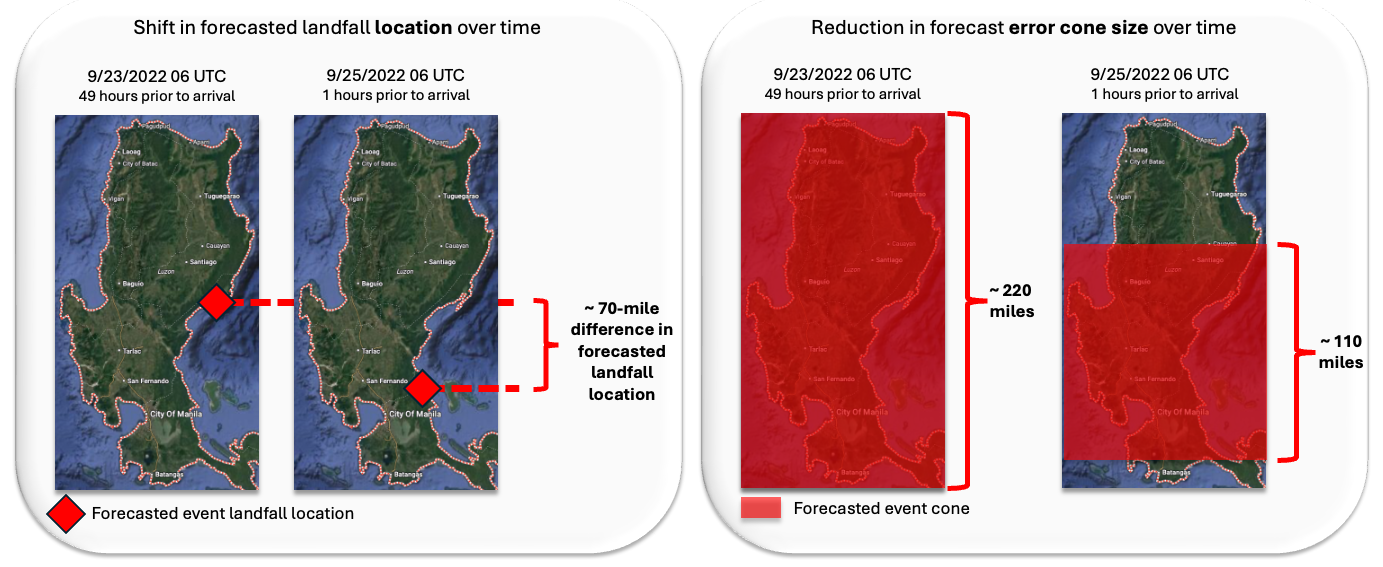}}
{The left panel shows the roughly 70-mile shift south in the forecasted event landfall location between early- and late-stage forecasts. The right panel illustrates the corresponding forecast cone narrowing from a width of roughly 220 miles to 110 miles. \label{fig:forecast_quantification}}
{Map data ©2025 \cite{GoogleMaps2025} and \cite{ecmwf_datasets}} 
\end{figure}

Two changes in the forecast reveal themselves as especially important for understanding the timing trade-off, as shown in Figure~\ref{fig:forecast_quantification}. The first is the correction in the storm's expected trajectory, which resulted in a ~70-mile shift in the predicted event landfall location between the 49-hour and 1-hour forecasts. The second is the reduction in variance, in which the diameter of the forecast cone, representing the range of possible event arrival zones, contracted from roughly 220 miles to 110 miles. It is precisely this quantifiable improvement in forecast precision that the model evaluates against the diminishing lead time for deployment.

\subsubsection{Demand estimation}

\begin{figure}[!htb]
\FIGURE
{\includegraphics[width=0.4\linewidth]{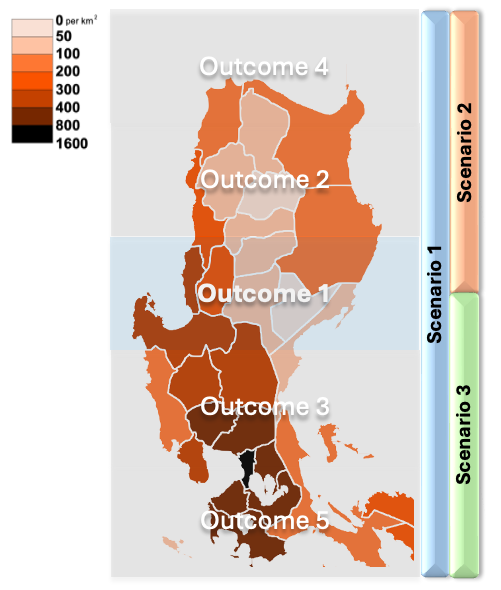}}
{Population density map of Luzon used as a proxy for the geographic distribution of demand.
\label{fig:demand_outcomes}}
{Source: \cite{WikipediaDemographicsPhilippines}}
\end{figure}

Humanitarian demand was estimated for 32 cities and towns located within the forecast cone. Demand is assumed to be proportional to the local population, which was derived from census data. Figure~\ref{fig:demand_outcomes} illustrates the population density and the discretization of the coastline into five potential event landfall zones.

The logistics network allows any demand location to serve as a candidate pre-staging site, which represents the use of local facilities as temporary distribution hubs. We estimate road travel times using a commercial routing application programming interface (API) that returns travel-time estimates between origin and destination locations. To reflect common relief commodities and to incorporate differing international lead times, we model two supply types, water and shelter, and source them from external hubs in the United States and Australia, respectively. Shelter provides temporary housing and protection, and potable water supports basic survival and public health during the first hours of response. We set the dispatch lead time to 24 hours from the United States and 12 hours from Australia, which represents the time from international dispatch to supply availability at the selected pre-staging sites.

\section{Results | Probabilistic optimization for dispatch timing}\label{sec:Results}

Table~\ref{tab:combined_results} integrates ($\mathcal{P}_{\!1}$) and post-processing results with the ($\mathcal{P}_{\!2}$) evaluation for water sourced from the United States. The results are presented as the expected closing time penalty and the expected unmet demand, which together form the objective value. The user-defined coefficients, $a_k$ and $b_k$, are  equal in the Table~\ref{tab:combined_results} results. 

\begin{table}[H]
  \centering
  \caption{Integrated ($\mathcal{P}_{\!1}$) Outputs and ($\mathcal{P}_{\!2}$) Objective Components for Typhoon Noru. The values in red in the (last four) solutions columns indicate suboptimality.}
  \label{tab:combined_results}
  \small
  \newcolumntype{C}[1]{>{\centering\arraybackslash}m{#1}}
  \begin{tabular}{C{2.2cm}|C{2.2cm}||C{2.2cm}|C{2.2cm}|C{2.2cm}|C{2.2cm}}
    
    \textbf{Dispatch Time} & \textbf{Expected Closing} & \textbf{Lateness} & \textbf{Unmet Demand} & \textbf{Time Penalty} & \textbf{Objective Value} \\
    (hrs prior) & (hrs) & (hrs) & (\%) & (hrs) & (composite) \\
    \hline
    49 & 106 & 0  & \textcolor{red}{29} & \textcolor{red}{106} & \textcolor{red}{135} \\
    37 & 98  & 0  & \textcolor{red}{28} & \textcolor{red}{98}  & \textcolor{red}{127} \\
    25 & 75  & 0  & \textbf{20} & \textbf{75} & \textbf{96} \\
    13 & 71  & \textcolor{red}{11} & \textcolor{red}{15} & \textcolor{red}{82}  & \textcolor{red}{97} \\
    1  & 66  & \textcolor{red}{23} & \textcolor{red}{9}  & \textcolor{red}{89}  & \textcolor{red}{98} \\
    
  \end{tabular}
\end{table}

\begin{figure}[H]
\FIGURE
{\includegraphics[width=\linewidth]{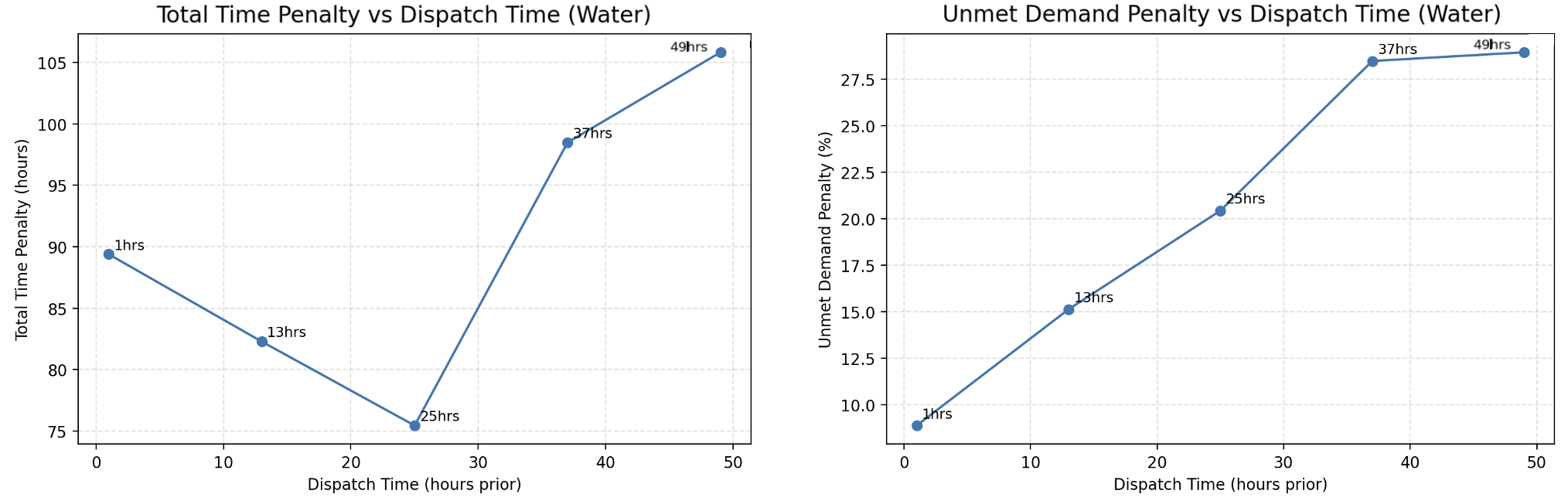}}
{($\mathcal{P}_{\!2}$) solution for water displaying the trade-off between the total time penalty and unmet demand.
\label{fig:water_results}}
{}
\end{figure}

The results show the trade-off between dispatching early under greater forecast uncertainty and waiting for improved information at the cost of reduced lead time. As the dispatch time shifts from 49 hours to 1 hour before landfall, expected closing time decreases by about 38 percent. However, dispatch times inside the 24-hour lead window incur a lateness penalty because supplies cannot reach the pre-staging sites before demand materializes. The objective reaches its minimum at a dispatch lead time of 25 hours; Figure~\ref{fig:water_results} demonstrates. 

\subsection{Varying weight parameters and sourcing locations}

The objective combines two terms expressed in different units, hours and percent unmet demand. This unit mismatch makes direct comparison subjective. For example, a stakeholder must decide whether a 15\% increase in the total time penalty justifies a 15\% reduction in unmet demand, or whether the same time increase justifies only a larger reduction in unmet demand. The model cannot resolve this judgment without stakeholder input.

\begin{figure}[H]
\FIGURE
{\includegraphics[width=0.7\linewidth]{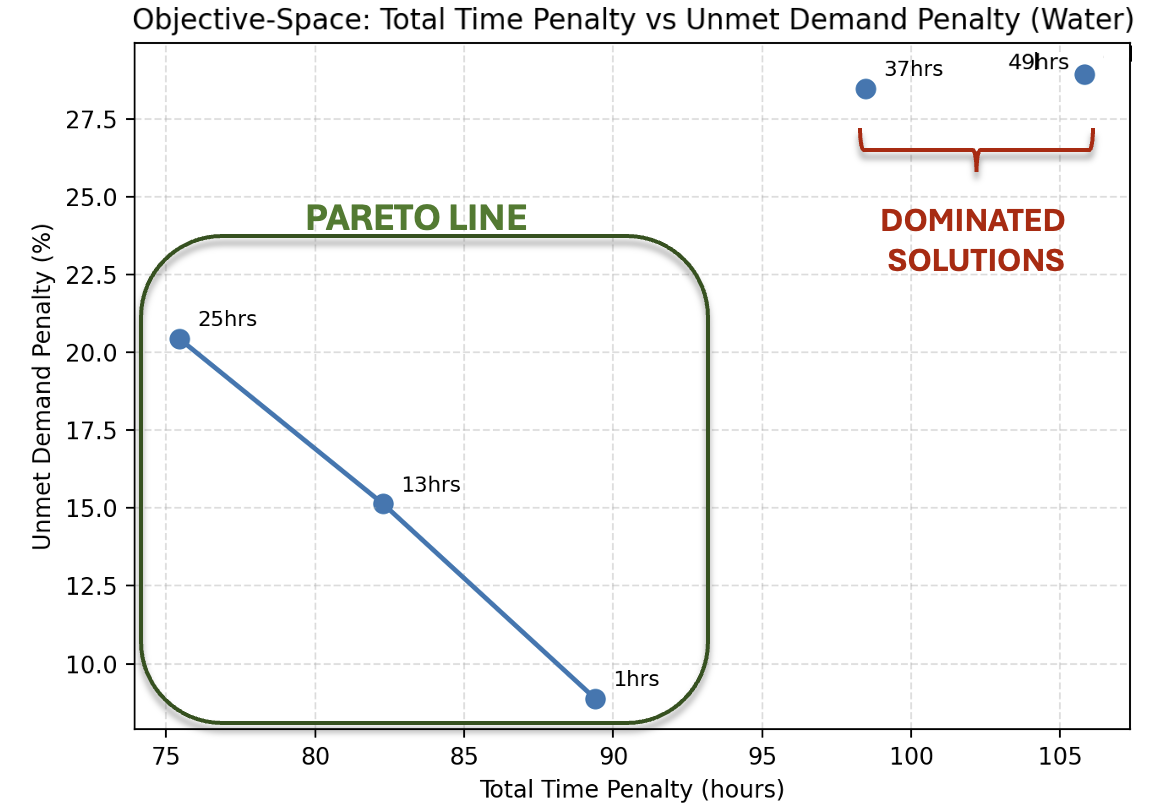}}
{Objective-space trade-off for water time penalty and unmet demand.
\label{fig:pareto_water}}
{}
\end{figure}

We recommend enumeration when the candidate set remains small or the models solve quickly because it makes trade-offs transparent without imposing an arbitrary conversion between hours and percent. For larger candidate sets, Pareto analysis would reveal dispatch trade-offs between  expected closing  time  and expected unmet demand. Figure~\ref{fig:pareto_water} plots these results across both objective function terms.

\begin{figure}[H]
\FIGURE
{\includegraphics[width=\linewidth]{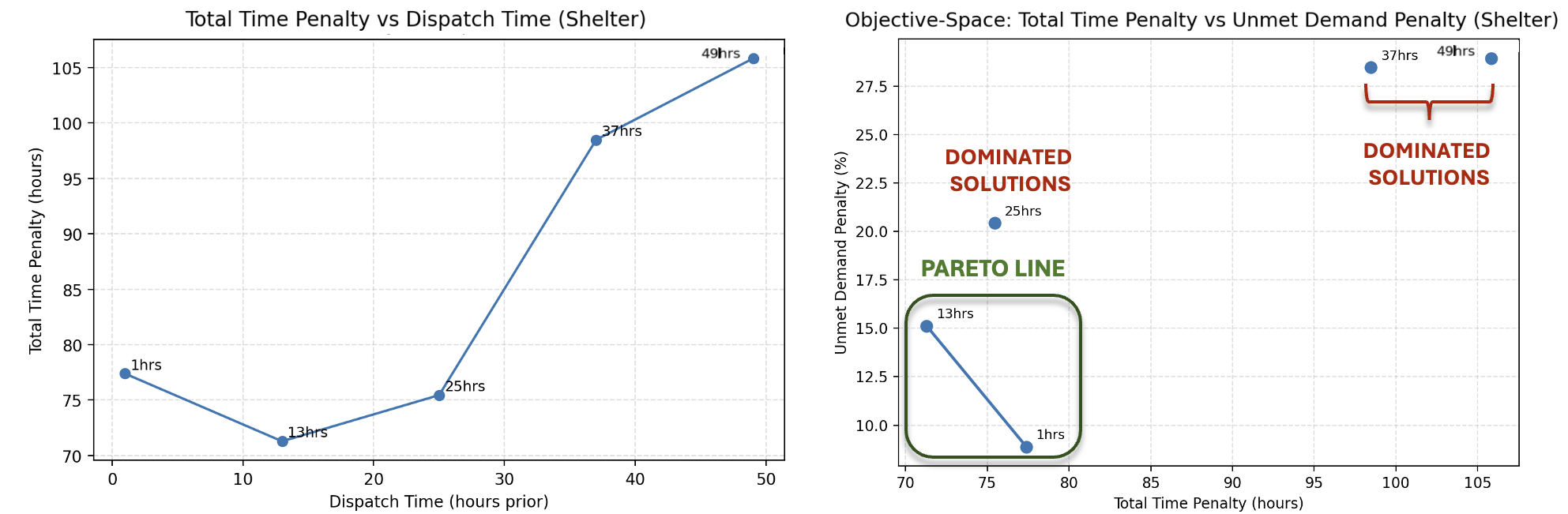}}
{Objective-space results for shelter with a 12-hour dispatch lead time showing the shift in the efficient frontier.
\label{fig:shelter_pareto}}
{}
\end{figure}

When shelter supplies originate from Australia, the shorter 12-hour lead time allows dispatch decisions to be deferred later in the forecast horizon. Indexing dispatch timing and performance metrics by sourcing location permits the model to capture the value of this added temporal flexibility and to evaluate how later commitment changes the trade-off between closing time penalties and unmet demand. Figure~\ref{fig:shelter_pareto} shows that a 13-hour dispatch yields both a lower total-time penalty and reduced unmet demand compared to earlier dispatches.

\subsection{Probabilistic optimization for dispatch timing | Case study summary}

Model ($\mathcal{P}_{\!1}$), post-processing, and model ($\mathcal{P}_{\!2}$) select a dispatch time that balances forecast uncertainty against lead-time loss. The transition from planning to execution occurs, in this case, 25 hours before event landfall. At that point, the organization can optionally shift to a higher-fidelity model for last-mile distribution, which routes deliveries from pre-stage sites to demand locations on the local road network.

\section{Robust optimization for last-mile distribution}
The preceding sections determine when humanitarian supplies should be dispatched and identify suitable staging locations at a coarse spatial scale under uncertain event forecasts. Once these upstream decisions are fixed, the focus shifts to local, i.e., last-mile, execution. In non-network-disrupted last-mile delivery, operational design choices, such as service zoning and vehicle routing, can materially affect efficiency, cost, and workload equity \citep{Manikas2016, Sinha2026, Carlsson2025}. Under disaster conditions, however, the dominant task shifts from simply minimizing closing time to whether supplies staged at otherwise-viable locations can reach affected populations through a damaged transportation network. In humanitarian logistics, the last-mile distribution step moves relief supplies from local distribution centers to beneficiaries, and connects routing and allocation decisions directly to infrastructure conditions under network damage \citep{Balcik2008LastMile, NoyanBalcikAtakan2015LastMile}.

\subsection{Network accessibility}
Although the selected pre-stage locations are assumed to remain operational under the event, the road network connecting those locations to demand is not explicitly modeled in the earlier phases. Last-mile delivery depends on a limited number of highways, arterial roads, and access routes that may be compromised by flooding, wind damage, or significant traffic. As a result, supplies can become isolated even when pre-staging sites themselves are located outside of the event trajectory and remain operational. Addressing this risk requires shifting attention from the survivability of locations to the accessibility of the network that links them to demand.

\subsection{Considering robustness at the last mile}

Unlike forecast trajectory uncertainty, last-mile infrastructure disruption is highly localized, environment-dependent, and difficult to characterize probabilistically. Road failures may result from flooding, debris, or structural damage that vary sharply over small spatial scales. Rather than attempting to assign failure probabilities to individual road segments, the primary operational concern becomes exposure to worst-case connectivity loss.

Robust optimization addresses this concern by seeking solutions that perform well under adverse realizations without requiring precise probability estimates. Foundational work by \cite{BenTalNemirovski1998} establishes tractable methods for handling uncertainty sets in linear and nonlinear programs. In humanitarian contexts, \cite{Najafi2013} develop a robust multi-objective model for earthquake response that ensures feasibility under demand uncertainty and network disruption, while \cite{Li2025RobustELCP} demonstrate improved resilience in emergency resource location–allocation decisions under multiple sources of uncertainty.

Recent extensions reduce reliance on full distributional assumptions in classical stochastic models. \citet{Yang2021ScenarioRobust} introduce scenario-robust prepositioning where planners specify distribution-free uncertainty sets in the form of ranges, rather than requiring exact probability distributions. \citet{Mahtab2022} propose a multi-objective, multi-period robust–stochastic humanitarian logistics model that accounts for uncertainty in demand, network reachability, and post-disaster inventory conditions, and includes equity considerations, demonstrated via a Bangladesh flood case study. \citet{DeMoorWagenaar2024} develop a multi-period robust and adaptive robust optimization model for humanitarian food aid operations under uncertain procurement prices, allowing some decisions to be made after uncertainty is revealed. \citet{Bertsimas2024} introduce coupled uncertainty sets to reduce the conservatism of robust and adaptive-robust optimization, and provide tight computable bounds. 

Despite this progress, robust methods in humanitarian logistics typically focus on facility location or post-event allocation, rather than on spatial diversification in last-mile routing. For this tutorial, we demonstrate a case study in which route supply concentration can isolate large pre-staged supplies from demand if part of the route becomes inaccessible. This setting provides a clear context for illustrating robust modeling choices that protect delivery performance under network-level disruptions. In this setting, the dominant risk is over-concentration. If a large share of supplies is staged in a single urban hub or routed through a small number of corridors, a local network failure can isolate a disproportionate fraction of inventory.

\subsection{Over-concentration of supplies: Outcome~3 illustration}

This risk is illustrated by Outcome~3 in the Typhoon Noru case study, which is the most likely landfall realization at the dispatch time recommended by ($\mathcal{P}_{\!2}$) (25 hours prior to arrival). Under this outcome, a large share of projected demand lies in Angeles, San Fernando, and San Miguel, north of Metro Manila. Within the solution provided by ($\mathcal{P}_{\!1}$), pre-staging a substantial portion of supplies in Manila is optimal within the solution provided by ($\mathcal{P}_{\!2}$): Manila minimizes expected travel time, provides efficient access to the affected region, and is located outside of the projected event trajectory at the time of dispatch.

However, access from Manila to the northern demand centers relies heavily on a single high-capacity highway corridor: the North Luzon Expressway running through Bocaue and Malolos. If this primary corridor is compromised by the event, supplies concentrated in Manila may become inaccessible via planned routes. While alternative routes may still allow delivery, detours can significantly increase closing time; if delivery is not feasible, demand may be left unmet. This risk is particularly acute for time-sensitive commodities such as potable water and shelter supplies, for which delayed delivery directly translates into preventable harm.

\section{Methodology | Robust optimization for last-mile distribution}

To mitigate the risk of over-concentration, we adopt a robust optimization approach based on an iterative adversarial loop \citep{Bienstock2010}. In this framework, a ``Blue'' decision-maker (the logistical planner) selects a dispatch plan, while a ``Red'' adversary identifies network vulnerabilities by seeking to maximize the consequences of Blue's plan through network disruptions. Critically, the Red model does not correspond to a specific storm forecast or scenario; rather, it systematically identifies the road segments whose failure would cause the greatest disruption given the current routing decisions.

The algorithm proceeds iteratively in a feedback loop: at each iteration, the Blue planner computes a routing plan, the Red adversary evaluates that plan to identify vulnerabilities, and the resulting penalties serve as inputs  into the next Blue solve:

\begin{enumerate}
    \item \textbf{Baseline optimization:} The Blue agent solves a mixed-integer program (MIP) to minimize delivery time, assuming a fully intact network. This produces a lowest-cost solution with no awareness of over concentration. This model is labeled ($\mathcal{R}^{-}$) and its solution is referred to as the {\it risk-unaware} solution.
    \item \textbf{Vulnerability assessment:} The Red agent analyzes the resulting traffic flow and identifies arcs with the highest vehicle concentration, i.e., the ``heavy links.''
    \item \textbf{Penalty injection:} The costs of these heavy links are increased in the network graph, representing the ``risk cost'' of over-reliance on a single road link.
    \item \textbf{Reoptimization:} The Blue agent re-solves the MIP on the penalized network. Because the primary highway is now more costly, the solver diverts a portion of the traffic to secondary roads or parallel highways. This model is labeled ($\mathcal{R}^{+}$) and its solution is referred to as the {\it robust solution}. 
\end{enumerate}

By iteratively adjusting the dispatch plan in response to these worst-case disruptions, the model gradually reduces reliance on any single location or corridor. The resulting solution does not attempt to equalize allocations, and it does not impose fixed upper bounds on link flows or on inventory staged at any site. Instead, diversification emerges endogenously as a response to identified vulnerabilities, producing a last-mile dispatch plan that is resilient to a broad class of severe, but possible, infrastructure failures.

\subsection{Robust route-based formulation}

\begin{figure}[htbp]
    \centering
    \includegraphics[width=0.95\linewidth]{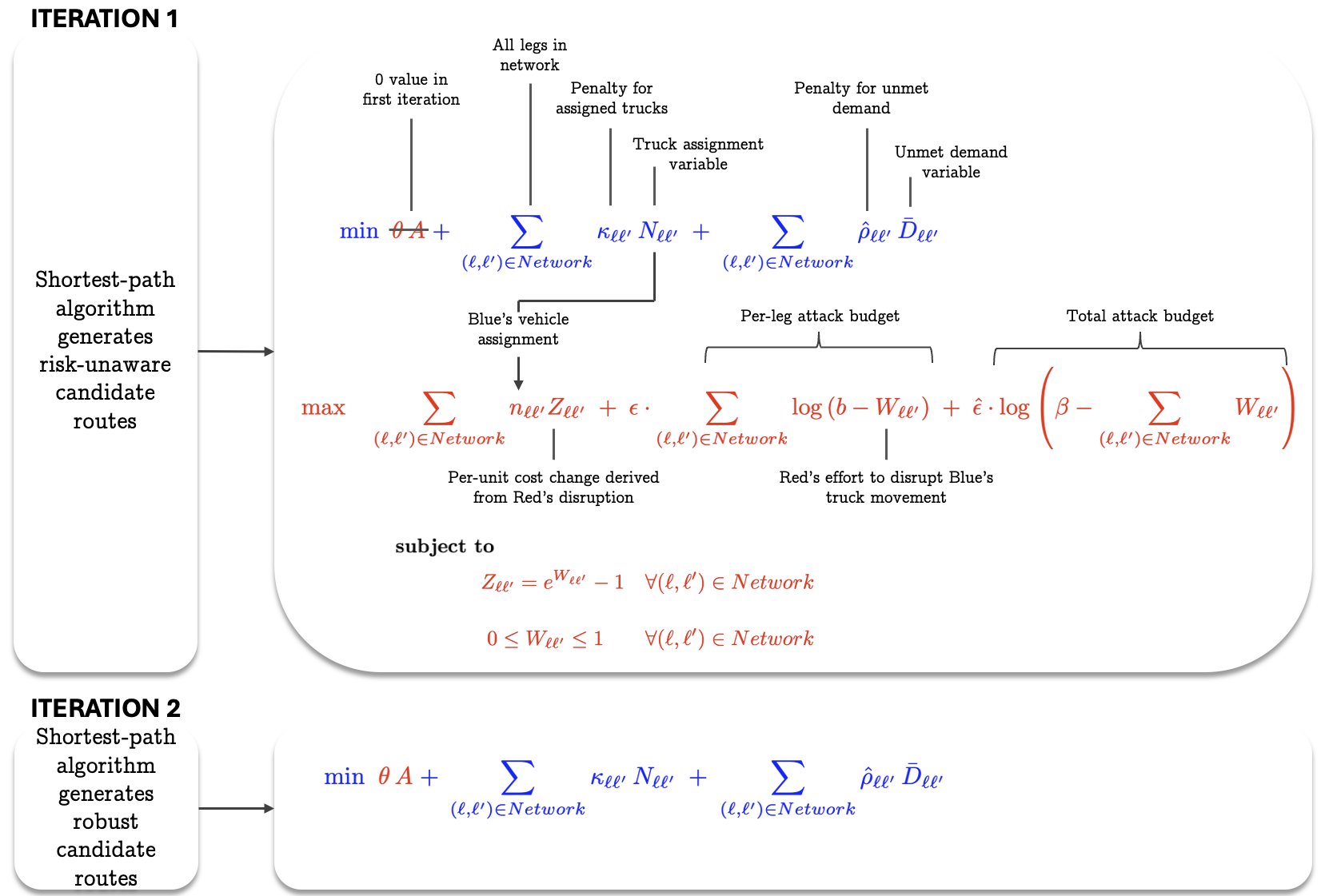}
    \caption{Schematic of the iterative Blue--Red robust-routing loop. Parameters: $\theta$ weights Blue’s exposure penalty; $\kappa_{\ell\ell'}$ is the nominal traversal cost of link $(\ell,\ell')$; $\hat{\rho}_{\ell\ell'}$ is the penalty applied to unmet demand; $\beta$ is Red’s total disruption budget; $b_{\ell\ell'}$ is Red’s per-link disruption budget; $\epsilon$ and $\hat{\epsilon}$ are log-barrier weights enforcing the per-link and total budgets, respectively. Red’s decision $W_{\ell\ell'}$ produces a per-unit cost increase $Z_{\ell\ell'}=e^{W_{\ell\ell'}}-1$, which is injected into the next Blue solve.}
    \label{fig:red_blue_schematic}
\end{figure}

Figure~\ref{fig:red_blue_schematic} summarizes a simplified version of the iterative robust last-mile formulation. The Blue model is more detailed than ($\mathcal{P}_{\!1}$) because ($\mathcal{P}_{\!1}$) operates at a regional level and represents movement with aggregate origin-destination travel times. Here, Blue instead enforces last-mile flow and assignment constraints on an explicit route-based road network, so each shipment follows a route composed of multiple road links, rather than a single origin-destination link.

Blue's decision-making occurs in two phases within each iteration. First, a shortest-path algorithm generates candidate routes for each origin-destination pair. When road links with many vehicles have been identified during previous vulnerability assessments, the Red model penalizes these arcs by adding a risk term to their traversal cost. If this penalty exceeds the cost of a longer detour, the algorithm discovers an alternative route. Second, Blue's MIP selects among all routes to minimize an objective that includes: (i) a risk term penalizing vehicle exposure on heavy links, (ii) vehicle usage costs, and (iii) unmet demand penalties. Although this formulation is different than the ($\mathcal{P}_{\!1}$) model, these penalties still favor rapid closing time. We present the full Blue and Red formulations in the appendix. However, in Figure~\ref{fig:red_blue_schematic}, we intentionally simplify and generalize the index notation to focus on the iterative structure of the algorithm.

\subsection{Red's nonlinear disruption}

Red’s disruption model is designed to capture the empirical observation that the operational consequences of targeted failures can grow nonlinearly when a system is highly concentrated. The variable $W_{\ell\ell'}$ represents Red’s effort (or severity) applied to link $(\ell,\ell')$, and the transformation $Z_{\ell\ell'} = e^{W_{\ell\ell'}}-1$ maps that effort into a per-unit cost increase that Blue experiences in the next iteration. This exponential mapping encodes the idea that incremental disruption can have \emph{out-of-proportion} impacts once a critical link is disabled. Exponential growth is not the only choice, but it provides a simple way to model compounding effects of over-reliance. This behavior is consistent with real-world supply-chain disruptions in which a single chokepoint failure produces cascading impacts, such as the 2021 \emph{Ever Given} incident in which the container vessel had run aground, blocking the Suez Canal for six days and causing severe global repercussions \citep{WAN2023106868}.

\subsection{Red budget limits and diminishing returns}
The log-barrier terms in Red’s objective play a complementary role by enforcing disruption budgets while reflecting diminishing marginal returns as Red approaches its limits. Intuitively, Red can increase disruption on individual links and also increase total disruption across the network, but doing so becomes progressively ``harder'' near the imposed limits: as Red tries to spend nearly all of its per-link or total budget, the model imposes a steep penalty that makes further concentration unattractive. This captures a practical reality in both adversarial and hazard-driven settings: extreme disruption is resource constrained, and additional effort near a limit often yields smaller incremental gains than earlier investments. Under this interpretation, Red’s ``effort'' variable can represent either intentional adversarial action or localized hazard intensity, which makes the same robust-routing loop applicable beyond the Philippines case study. As with the nonlinear mapping from Red’s effort to the disruption Blue perceives, the log barrier function is only one reasonable choice. In practice, the user can select any penalty function that matches how quickly disruptions should intensify and how sharply budget limits should constrain Red as those limits are approached.

\subsection{Determining $\theta$ in the risk term}

In this setting, the number of delivery trucks is approximately proportional to the quantity of supplies moved, so truck flow provides a natural proxy for exposure. However, Red’s first objective term is defined in terms of aggregate flow on links and is therefore not restricted to truck traffic; it can represent any exposure measure the user chooses (e.g., tonnage, passengers, value, or mission-critical throughput).

The risk term in Blue's objective takes the form $\theta \cdot A$, where $A$ represents the maximum vehicle exposure across all attack scenarios identified by Red, and $\theta$ is a scalar parameter that controls the relative importance of risk avoidance. This feature is similar to the risk aversion term in the Markovian Markowitz Mean Variance method \citep{Markowitz1952PortfolioSelection} in which a risk-aversion coefficient weights the variance penalty relative to expected return, thereby controlling the trade-off between portfolio performance under nominal conditions and protection against adverse outcomes. 

The parameter $\theta$ influences the algorithm at two distinct points. First, during route generation, the shortest-path algorithm uses $\theta$ to increase the traversal cost of each heavy link. If $\theta$ exceeds the marginal cost of the next-shortest path outside the current route pool, the algorithm will discover and add an alternative route that avoids the concentrated link. In this case, new routes enter the pool, expanding Blue's options for subsequent iterations.

Second, during route selection, the MIP balances the risk term $\theta \cdot A$ against the vehicle usage and unmet demand penalties. Even if $\theta$ is not large enough to generate new routes, the optimizer may still diversify within the existing route pool. This occurs when $\theta \cdot A$ exceeds the marginal cost of reassigning vehicles to less-exposed routes or shifting demand fulfillment from origins with heavy-link exposure to origins with cleaner routes. In this way, the MIP can reduce risk by redistributing which origins serve which destinations, even when each origin-destination pair has only a single route available.

Thus, $\theta$ serves a dual role: it determines whether new attack-avoiding routes are generated, and it controls how extensively the optimizer diversifies among routes already in the pool. Setting $\theta$ too small results in minimal or no diversification; setting $\theta$ too large may sacrifice efficiency by overweighting risk avoidance relative to vehicle costs and demand fulfillment. In practice, we recommend initializing $\theta$ to a value larger than the longest route length in the network and inspecting the solution after the first iteration. If new routes are generated but the optimizer does not use them, the vehicle usage and unmet demand penalties in the MIP may dominate the risk term; increasing $\theta$ will strengthen the incentive to diversify. This parameterization provides the user with flexibility: rather than imposing mandatory diversification, $\theta$ acts as a knob that the user can adjust to explore solutions along the diversification-efficiency frontier. Without this knob, the user would have no ability to control the tradeoff between resilience and operational cost.

\subsection{A note on alternative methods to diversify the Blue MIP}
Alternatively, a planner could attempt to prevent over-concentration by adding explicit flow or warehouse capacities, such as knapsack-style constraints that limit flow on each road link or limit inventory at each staging location. However, this approach requires the planner to specify upper bounds that are difficult to justify because the appropriate capacity depends on network topology and on how disruption consequences propagate through connectivity. For example, if the upper bounds are aggressive, i.e., too small, the optimization problem may be rendered infeasible, or nearly so; and, in the latter case, the solution may exhibit undesirable artifacts. 

Rather than imposing fixed upper bounds, the vulnerability assessment identifies which concentrations are problematic given current routing decisions, and the penalty injection mechanism raises the cost of those specific road links. This process adapts to network topology automatically: road links on multiple near-optimal routes receive higher penalties as traffic accumulates, while routes with alternatives see traffic diverted before penalties grow large. This process ensures that diversification emerges only when beneficial alternatives exist; if an alternative route is too costly, it is not generated, and flow remains concentrated as a deliberate outcome rather than as an artifact of arbitrary bounds. The result is diversification that is topology-aware and cost-justified, without requiring the planner to specify link-level or site-level limits. 


\section{Results | Robust optimization for last-mile distribution}
\begin{figure}[htbp]
    \centering
    \includegraphics[width=1.0\textwidth]{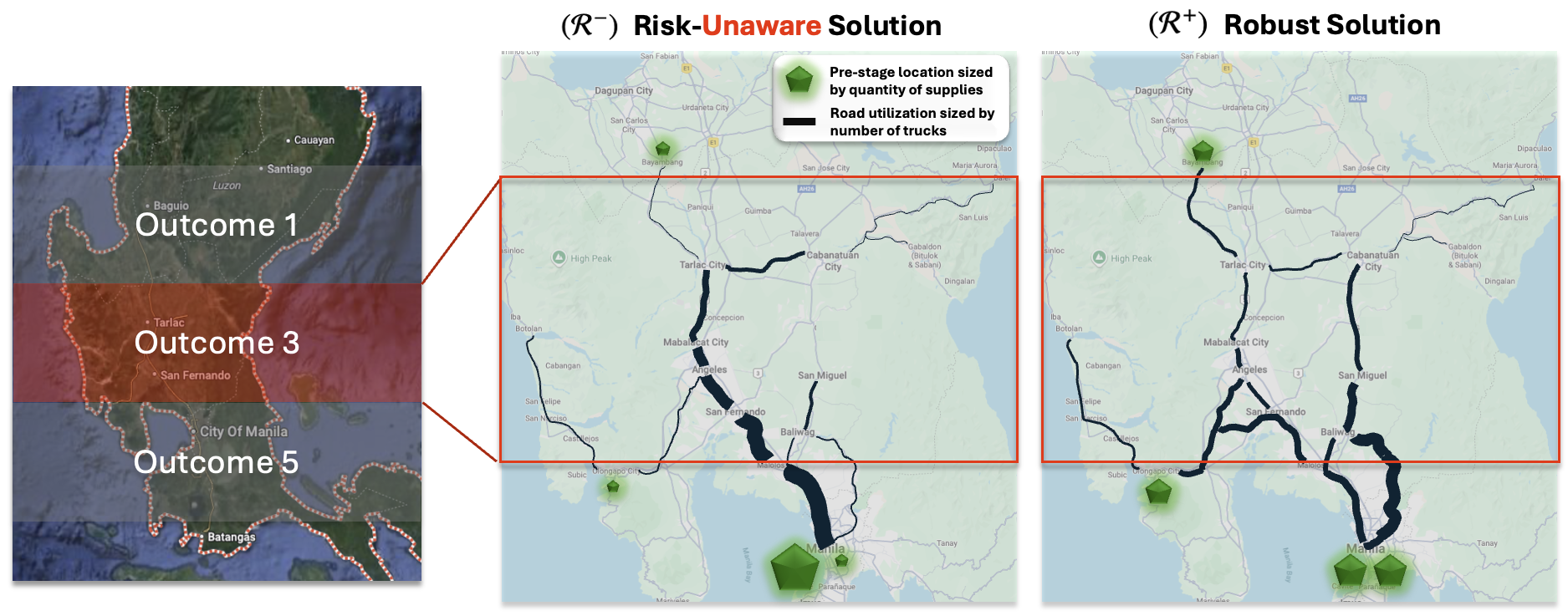}
    \caption{The solution provided by ($\mathcal{R}^{-}$) (left) channels the majority of trucks through the central corridor (Manila$\to$Bocaue$\to$Malolos$\to$San Fernando$\to$Angeles). The solution provided by ($\mathcal{R}^{+}$) (right) reduces traffic on this primary artery and activates two alternative routes: a western route through Floridablanca and an eastern route via Quezon City, Norzagaray, Baliwag, and San Miguel.}
    \label{fig:spatial_map}
\end{figure}

We evaluate the robust routing method under Outcome~3, which represents the most likely landfall realization at the dispatch time selected by ($\mathcal{P}_{\!2}$). In this analysis, we model last-mile delivery with a single vehicle mode, i.e., truck, for both commodities and map aggregate link flows to assess diversification. We run a sensitivity analysis and select a sufficiently large value of $\theta$ to initiate new routing options and flow diversification by the second iteration. Figure~\ref{fig:spatial_map} contrasts vehicle concentration in the solution provided by ($\mathcal{R}^{-}$) (left) and ($\mathcal{R}^{+}$) (right). Traffic volumes report the total number of truck trips on each road segment.

In the solution provided by ($\mathcal{R}^{-}$), the solver routes the majority of traffic through the central corridor, as this represents the shortest path to Angeles and San Fernando. The iterative vulnerability assessment identifies this concentration and activates two alternative corridors:

\begin{itemize}
    \item \textbf{Western Corridor}: Traffic increases from 1 to 25 trucks. This route provides an alternative route to Angeles and San Fernando that bypasses the congested Malolos$\to$San Fernando segment.
    \item \textbf{Eastern Corridor}: Traffic increases from 18 to 61 trucks. This corridor provides independent access to San Miguel and Cabanatuan, with onward connections to Tarlac City.
\end{itemize}

\subsection{Implications for pre-staging strategy}

The solution provided by ($\mathcal{R}^{+}$) increased utilization of the western and eastern corridors and suggests that distributing pre-staged inventory across multiple supply hubs, specifically, increasing allocations at Olongapo City (western) and Quezon City (eastern), would further improve system resilience. This diversification of both \textit{routes} and \textit{pre-staging locations} provides resilience against infrastructure failures.

\subsection{Verifying the redistribution}

\begin{figure}[htbp]
    \centering
    \includegraphics[width=1.0\textwidth]{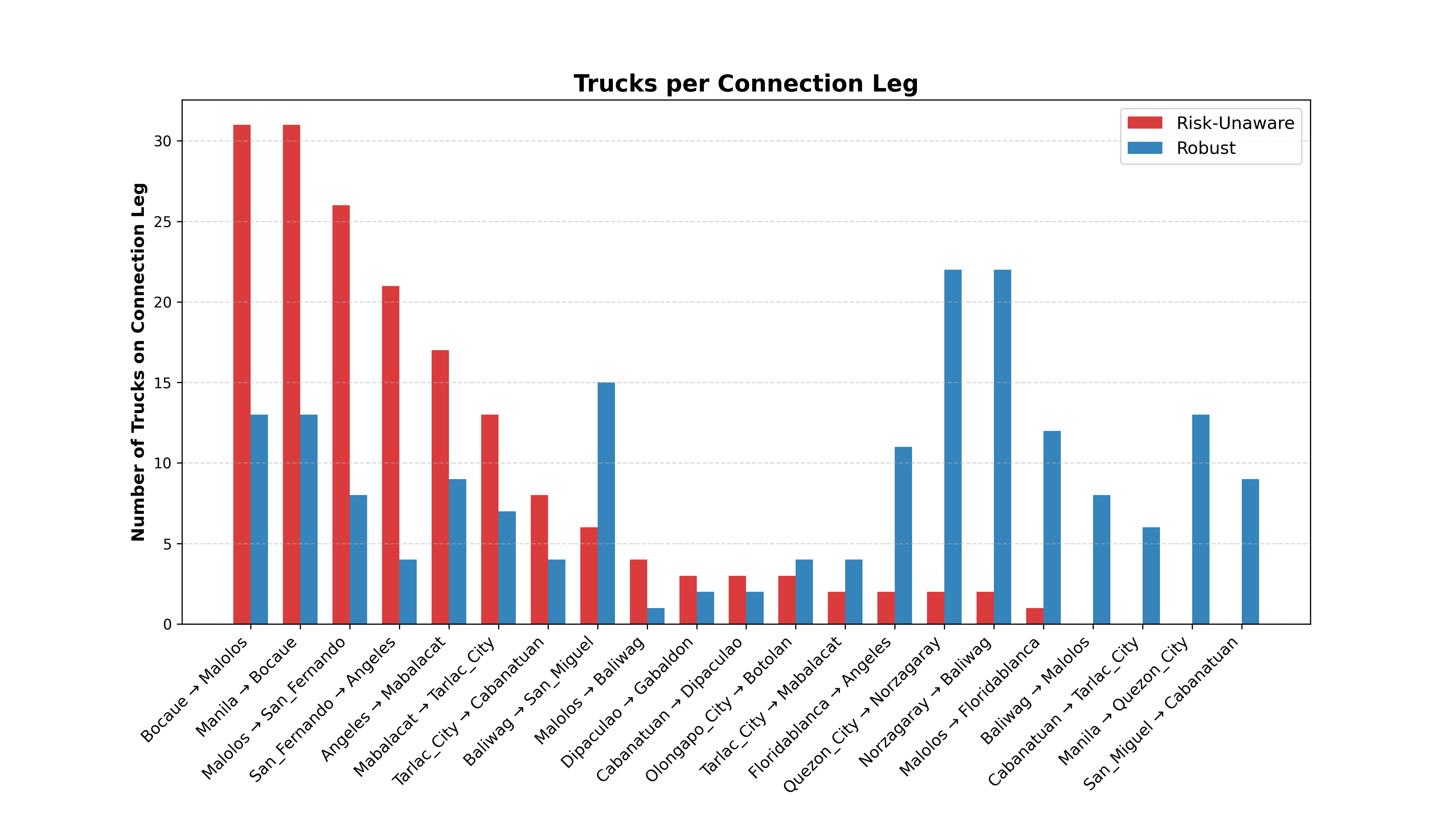}
    \caption{The solution provided by ($\mathcal{R}^{-}$) shows high peaks on the central corridor, with 31 trucks each on Manila$\to$Bocaue and Bocaue$\to$Malolos. After three iterations, the solution provided by ($\mathcal{R}^{+}$) reduces these peaks and redistributes load to the Norzagaray$\to$Baliwag segment (22 trucks) and the Floridablanca links. The solution provided by ($\mathcal{R}^{+}$) still concentrates some flow on the eastern corridor, but it reduces the maximum load and spreads traffic across multiple independent routes.}
    \label{fig:traffic_dist}
\end{figure}

To verify that the algorithm achieves diversification, rather than shifting the bottleneck elsewhere, we analyze the load distribution across all network links. Figure~\ref{fig:traffic_dist} plots the number of trucks on each road segment for both solutions and shows a flatter distribution under the solution provided by ($\mathcal{R}^{+}$). The solution curve provided by ($\mathcal{R}^{+}$) exhibits fewer extreme peaks, which indicates that traffic spreads across more links instead of concentrating on a small set of corridors. Table~\ref{tab:risk_metrics} summarizes the performance metrics in which we compare each solution's own top links to provide a fair assessment of concentration risk. 


\begin{table}[htbp]
\centering
\caption{Supply Concentration Metrics Provided by Models ($\mathcal{R}^{-}$) and ($\mathcal{R}^{+}$)$^{\dagger}$}
\label{tab:risk_metrics}
\begin{tabular}{lrrr}
\toprule
\textbf{Metric} & \textbf{($\mathcal{R}^{-}$)} & \textbf{($\mathcal{R}^{+}$}) & \textbf{Change} \\
\midrule
Peak single-link load (trucks) & 31 & 22 & \textbf{\textcolor{darkgreen}{$-$29\%}} \\
Top-3 links aggregate (trucks) & 88 & 59 & \textbf{\textcolor{darkgreen}{$-$33\%}} \\
Top-5 links aggregate (trucks) & 126 & 85 & \textbf{\textcolor{darkgreen}{$-$33\%}} \\
Central corridor total (trucks) & 109 & 38 & \textbf{\textcolor{darkgreen}{$-$65\%}} \\
\bottomrule
\end{tabular}
\vspace{3pt} \\
\parbox{\linewidth}{\footnotesize Note: $^{\dagger}$run for three iterations}
\end{table}

Across results provided by ($\mathcal{R}^{-}$) and ($\mathcal{R}^{+}$), the peak single-link load decreases 29\% (31 to 22 trucks). When comparing the aggregate load on each solution's three most-utilized links, the reduction is 33\% (from 88 to 59 trucks). The central corridor experiences the largest change, with a 65\% reduction in total traffic as the algorithm shifts load to the western and eastern alternatives (109 to 38 trucks).

Notably, the solution provided by ($\mathcal{R}^{+}$) does create new concentrations, the Norzagaray$\to$Baliwag segment now carries 22 trucks. However, this segment lies on an independent route from the central corridor. A failure on the central highway would not affect the eastern corridor, and vice versa. The goal of the robust methodology is not to eliminate all concentration, but to ensure that no single infrastructure failure can strand the entire relief operation.

\subsection{Delivery efficiency}

\begin{figure}[htbp]
    \centering
    \includegraphics[width=0.8\textwidth]{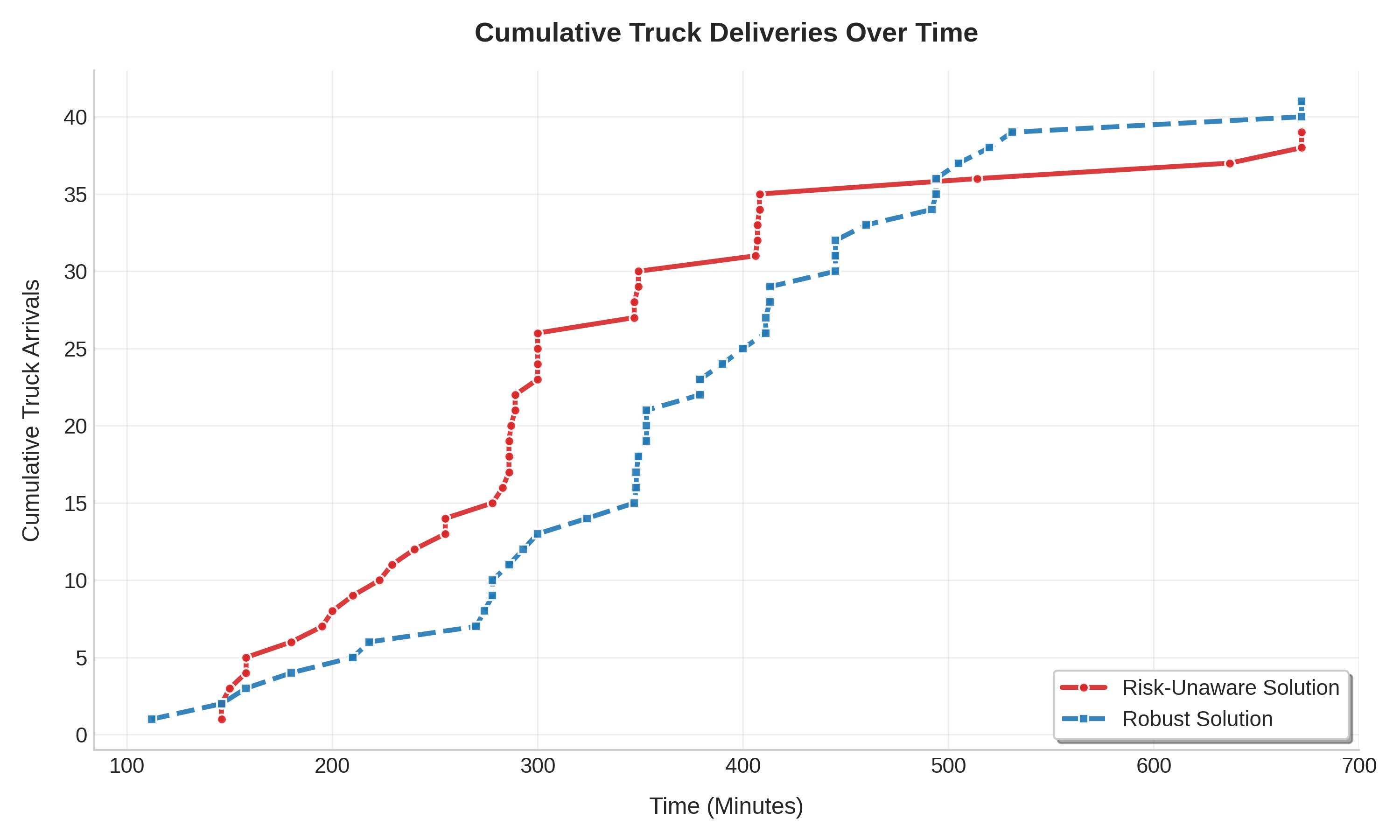}
    \caption{The solution provided by ($\mathcal{R}^{-}$) (red) achieves a steeper initial delivery curve. The solution provided by ($\mathcal{R}^{+}$) (Blue) lags at the beginning of the dispatch but converges by approximately 400 minutes into the operation. Both solutions complete 90\% of deliveries by roughly 350--400 minutes. The modest divergence indicates that the network offers sufficient redundancy to support spatially diverse routing without imposing a large delay penalty.}
    \label{fig:cumulative_deliveries}
\end{figure}

Figure~\ref{fig:cumulative_deliveries} compares the cumulative delivery timeline of both solutions.
The solution provided by ($\mathcal{R}^{-}$) delivers supplies faster at the beginning of the scenario. However, the gap narrows over time, and both solutions complete approximately 90\% of deliveries within a similar window (350--400 minutes). The final deliveries in both cases correspond to remote destinations (such as Gabaldon) that have limited routing alternatives, regardless of the optimization strategy.

The relatively small divergence between the curves suggests that the cost of diversification in this network, while subjective, is modest. The road network offers enough redundancy that the ($\mathcal{R}^{+}$) model can identify safer, spatially diverse routes without incurring a prohibitive time penalty. In operational terms, planners can achieve meaningful risk reduction while accepting only a minor increase in average delivery time. We note that different metrics can serve as proxies for robustness, including, for example, route diversity measured by statistical dispersion of link flows, maximum flow concentration on a single link, or delivery loss under a single-link failure, defined as the reduction in total delivered demand when the most heavily used segments become unavailable. We do not claim that any particular metric outcome is ``sufficient''; the acceptable tradeoff between diversification and efficiency depends on the operational context and must be evaluated by the user. 

Concluding the robust optimization methodology, we now step back to a broader tutorial question: how the uncertainty structure and available data should guide the choice between a probabilistic and robust   model.

\section{Discussion | Modeling choice that leverages network and data structure}

\begin{figure}[t]
\FIGURE
{\includegraphics[width=\linewidth]{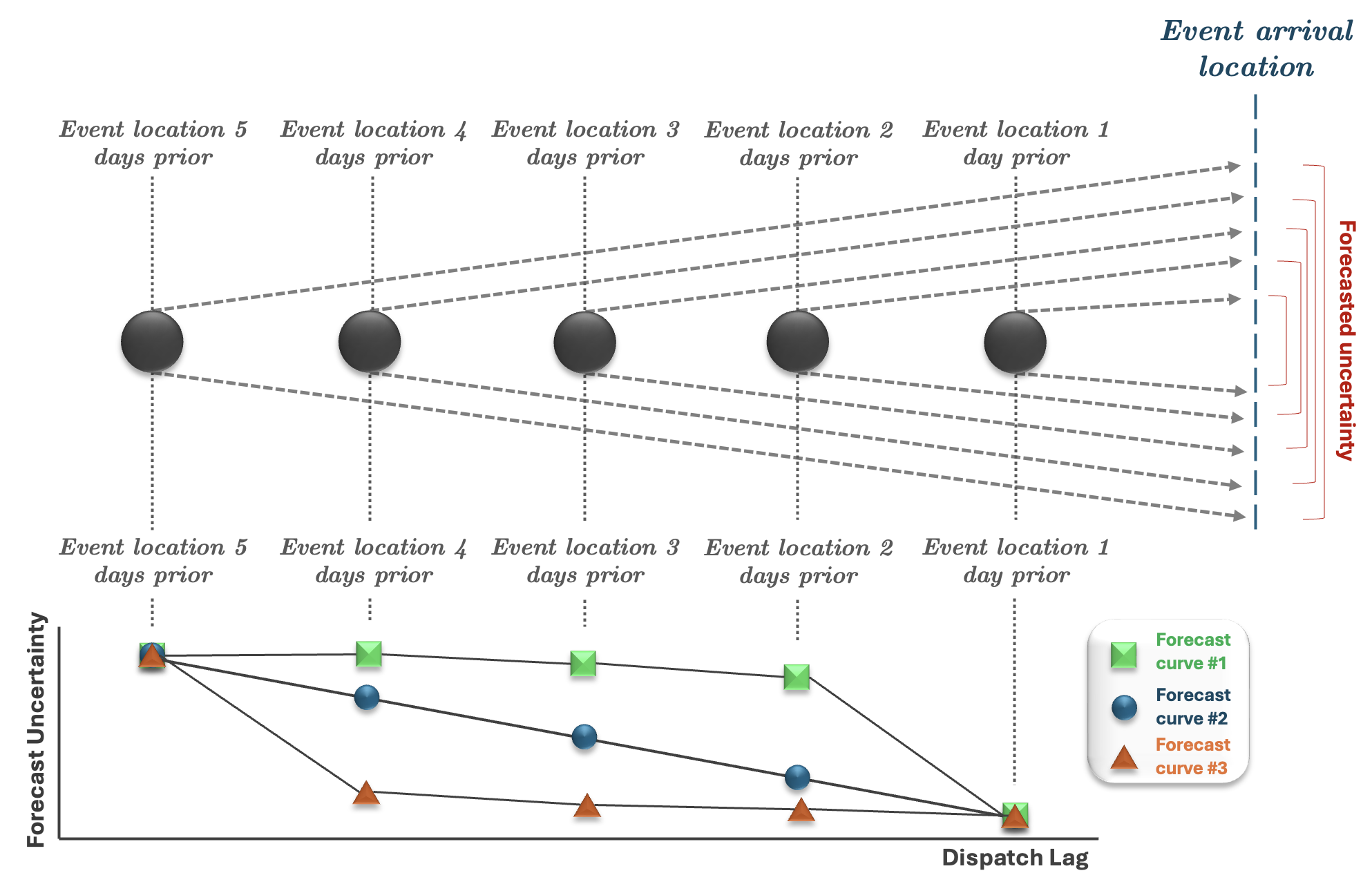}} 
{Generalized forecast-accuracy curves illustrating the trade space for dispatch timing. Each curve shows a potential shape of how forecast error evolves as an event nears the landfall area.
\label{fig:problemstatement}}
{Curves represent conceptual trends, not empirical data, used to illustrate timing trade-offs in anticipatory logistics decisions.}
\end{figure}

The ($\mathcal{P}_{\!1}$) and ($\mathcal{P}_{\!2}$) models and methodology demonstrate that an important  opportunity in anticipatory logistics lies in recognizing that the shape of the forecast-accuracy curve, rather than a single forecast snapshot, guides dispatch timing. Forecast accuracy does not improve uniformly across different global regions. As \citet{Ginis2022} document, western Pacific 48-hour track errors can be nearly twice those in the Atlantic; therefore, the value of waiting and the risk of acting early differ substantially by region. Figure~\ref{fig:problemstatement} illustrates this concept using generalized, conceptual forecast curves that represent how uncertainty evolves as an event approaches. Some curves tighten only shortly before landfall and others stabilize early. Each pattern defines a different trade space between early and delayed mobilization.

The two-phase framework leverages the specific timing of these improvements where ($\mathcal{P}_{\!1}$) provides a transparent pre-staging model that identifies spatial allocations under the current forecast. The ($\mathcal{P}_{\!2}$) model and ($\mathcal{P}_{\!1}$) post-processing transforms the evolving forecast curve into time-indexed expected impacts on closing time and unmet demand.

\subsection{Distributional topology and the limits of expectation}
\begin{figure}[t]
\FIGURE
{\includegraphics[width=\linewidth]{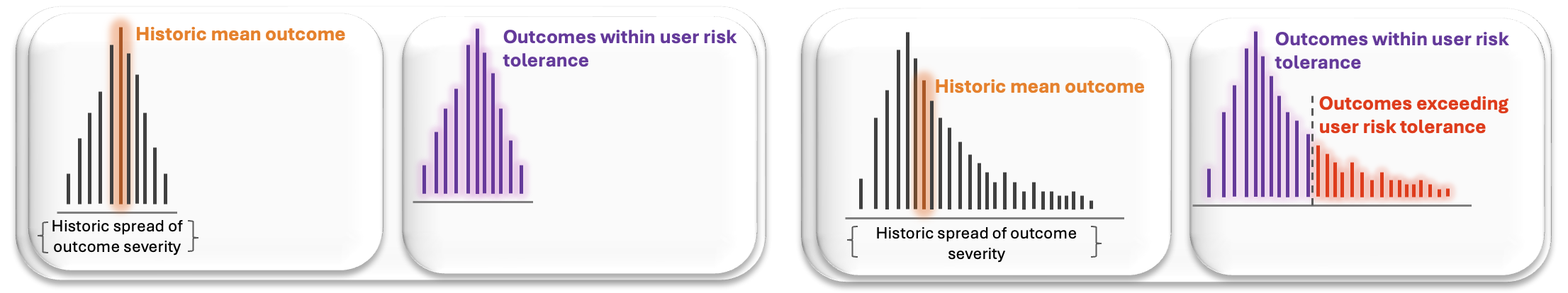}}
{Comparison of risk exposure under varying distribution shapes. The left panels show distributions whose outcomes remain within operational tolerance. The right panels show distributions where low-probability tails exceed system limits.
\label{fig:distributions}}
{}
\end{figure}

The probabilistic methodology evaluates performance in expectation, but its usefulness depends on the width and shape of the underlying outcome distribution. Distribution width captures the range of plausible realizations, while shape determines how much probability lies in the tails versus near the central outcomes. Figure~\ref{fig:distributions} illustrates this distinction. Some distributions exhibit limited spread, while others are wide or heavy-tailed. Width alone does not determine operational risk because risk depends on how realizations translate into closing time and unmet demand.

In models ($\mathcal{P}_{\!1}$) and ($\mathcal{P}_{\!2}$), uncertainty stems from the forecast (e.g., storm track or intensity). As this uncertainty propagates through the logistics system, it induces a distribution of operational outcomes, such as closing time, across landfall outcomes and forecast scenarios. This distribution need not be symmetric or tightly concentrated: a small number of realizations can interact with network topology and operational constraints to produce disproportionately large delays, yielding skewed or heavy-tailed outcomes. For instance, closing time may remain near its mean for most outcomes, yet a few low-probability outcomes can cause a substantial increase in closing time. Optimizing expected performance can therefore understate operational exposure, particularly because reliance on the mean implicitly assumes repeated exposure to uncertainty. In humanitarian logistics, that assumption may not hold: large-scale disasters are infrequent, planners observe only a small number of realizations, and performance may not ``average out'' within a season. Consequently, a single realization can govern outcomes, and that realization may fall in the tail rather than near the mean. Even when tail events carry low probability, their consequences may be severe enough to dominate decisions, in which case a robust approach that protects against adverse tail outcomes may be more appropriate. The relevant distributional shape is therefore system-specific and arises jointly from the forecast uncertainty and the network structure that maps realizations into response performance.

Relying on metrics of expectation becomes more acute when the distribution itself is difficult to specify. In many last-mile contexts, planners cannot reliably characterize uncertainty using probability distributions. For example, estimating post-storm road network failure requires detailed knowledge of localized storm intensity, rainfall accumulation over time, terrain and drainage, road construction standards, maintenance conditions, traffic incidents, debris patterns, and cascading effects from nearby failures. These factors interact at spatial and temporal resolutions that exceed available data and modeling fidelity.
Because the operational consequences of tail events can be so severe, it is prudent to assess solution performance under adverse, but possible, realizations of uncertainty rather than relying solely on expected values. Many practitioners turn to robust optimization models for this purpose. 

\subsection{The price of robustness}
The cost incurred when moving from a solution provided by ($\mathcal{R}^{-}$) to a more resilient alternative is commonly referred to as the \emph{price of robustness} \citep{BertsimasSim2004PriceRobustness}. In optimization problems, it is the increase in the objective value of the underlying MIP when the planner enforces robustness relative to the ($\mathcal{R}^{-}$) objective value. In the case study, Figure~\ref{fig:cumulative_deliveries} illustrates this trade-off by showing a modest loss in early delivery speed under the solution provided by ($\mathcal{R}^{+}$), which represents the closing time objective term. While some critics argue that robust methods are overly conservative, the price of robustness is not intrinsic to the methodology; it depends on network topology and data. Some networks admit optimal or near-optimal diversification at low cost, while others impose a substantial trade-off.

In networked logistics systems, this price is often lower than expected. Many transportation networks exhibit significant combinatorial flexibility: multiple routes, staging patterns, and flow allocations achieve identical or near-identical objective values under nominal conditions. As a result, planners can frequently introduce meaningful diversification with only a small degradation in closing time or total cost.

\begin{figure}[!htb]
\FIGURE
{\includegraphics[width=\linewidth]{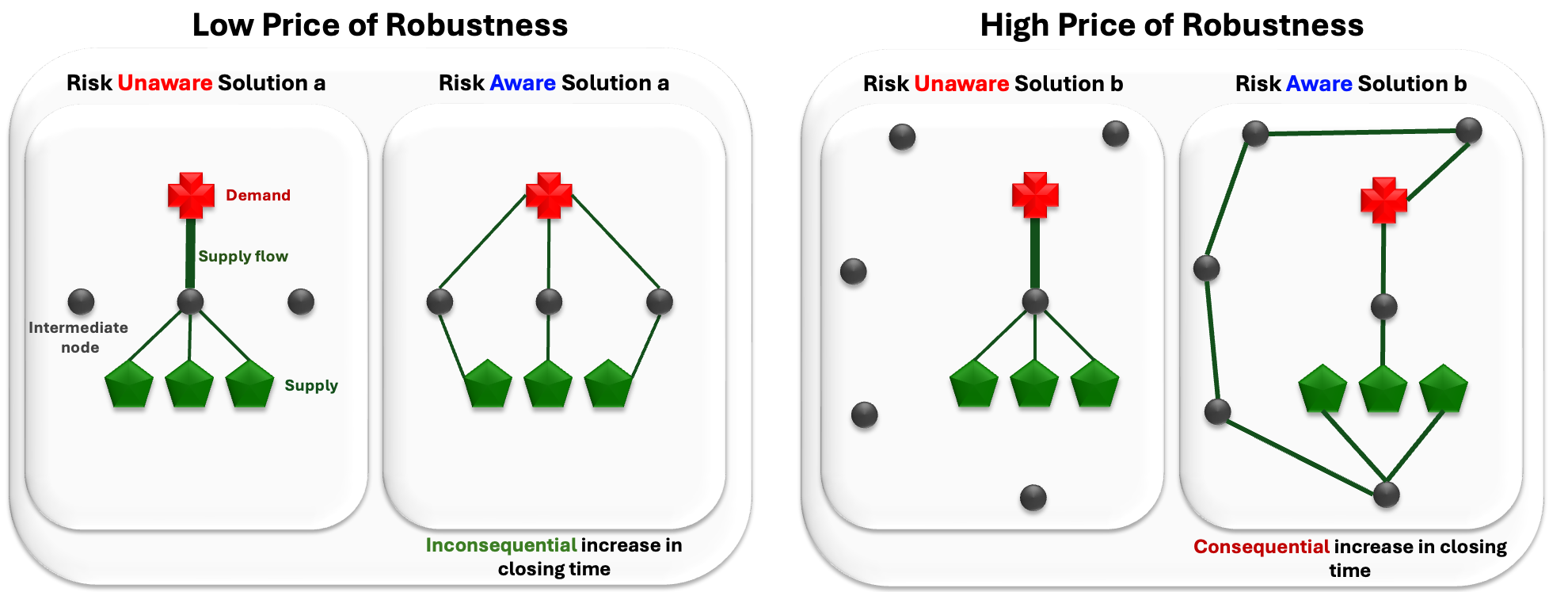}}
{ In each panel, the risk-unaware solution concentrates flow through the shortest, highly shared corridor, while the risk-aware solution deliberately spreads flow across multiple routes to reduce over-reliance on any single link. Panel (a) depicts a network with readily available, near-parallel alternatives, so diversification can be achieved with only a negligible increase in closing time. Panel (b) depicts a network with a single alternative route, so avoiding the bottleneck forces a long detour, resulting in a large increase in closing time (i.e., a high price of robustness).
\label{fig:price_robustness}}
{}
\end{figure}

Figure~\ref{fig:price_robustness} illustrates this concept, where a thicker supply flow line indicates a larger magnitude of flow. In Panel~(a), the solution provided by ($\mathcal{R}^{-}$) concentrates flow on the single shortest corridor, whereas the solution provided by ($\mathcal{R}^{+}$) diverts a portion of traffic to parallel routes-- increasing resilience to localized failures. The resulting increase in closing time is modest, indicating a low price of robustness. In Panel~(b), the network offers fewer viable alternatives. Introducing diversification requires routing flow along a long circumferential route, which substantially increases delivery time. In this case, the price of robustness is high. The methodology introduced in this tutorial exploits the first case. 

Because the robust methodology introduced in this tutorial relies on first-order methods, the resulting iterations remain computationally light and can be applied after any network optimization model. For this reason, a central recommendation of this tutorial is that practitioners  always evaluate the price of robustness. If a modest increase in closing time yields a disproportionate reduction in exposure to network failure, adopting the solution provided by ($\mathcal{R}^{+}$) is operationally justified. 

\clearpage
\section*{Conclusion}

This tutorial discusses probabilistic and robust optimization to address the lifecycle of humanitarian aid dispatch. It proposes a two-phase probabilistic model for determining the optimal time of dispatch and a robust model for securing last-mile distribution.

\subsection*{Probabilistic approach: optimizing dispatch timing}
The first section addresses the time-of-dispatch problem using a probabilistic lens, designed for strategic phases where uncertainty is well-characterized by probabilistic data such as weather forecasts. The methodology employs a two-phase framework:
\begin{itemize}
    \item Model ($\mathcal{P}_{\!1}$) addresses the spatial challenge by optimizing pre-staging locations based on the most likely storm trajectory.
    \item Model ($\mathcal{P}_{\!2}$) addresses the temporal challenge by evaluating the fixed plan against evolving forecast scenarios.
\end{itemize}
This approach identifies the specific time when the cost of lost lead time is outweighed by the value of improved forecast accuracy.

\subsection*{Robust approach: resilience in the last-mile delivery}
The second section addresses last-mile distribution using robust optimization. Once supplies are dispatched, the primary risk shifts from forecast uncertainty to difficult-to-predict infrastructure failures. Relying on standard cost minimization  often leads to over-concentration on single supply corridors.
To mitigate this, the tutorial introduces an iterative adversarial loop. A ``Blue'' agent plans efficient routes, while a ``Red'' agent penalizes bottlenecks to simulate disruption. This process forces the model to potentially, but not necessarily, accept longer delivery times in exchange for diversification, ensuring that no single road failure can significantly disrupt the relief effort.

\subsection*{Modeling choice}
Practitioners should utilize probabilistic optimization when the decision-maker trusts the probability distribution and accepts the risk inherent in the variability of outcomes. In this context, the deviation from the expected value is acceptable. Conversely, robust optimization is required when the tails of the distribution represent catastrophic failure or when the uncertainty structure cannot be reliably characterized, necessitating a plan that considers worst-case disruptions.
Selecting the correct modeling paradigm is not a philosophical exercise; it is a technical alignment of available data and unmeasured structure with the specific expectation that the decision-maker requires. This alignment transforms the organizational posture: it moves the firm from suffering unanticipated failures to anticipating variability and managing performance trade-offs before execution.

%
\clearpage

{\begin{APPENDIX}{}

\section{\large{Phase 1 Probabilistic Formulation, ($\mathcal{P}_{\!1}$)}}
\label{app:phase1}
\bigskip

{\fontsize{11}{12.5}\selectfont
\noindent\begin{tabular}{l l}
		\multicolumn{2}{l}{\bf \underline{Sets}} \\
		{\bf Symbol} & {\bf Definition} \\
		\hline
		$\mathcal{I}$ & potential pre-stage locations \\
		$\mathcal{J}$ & anticipated demand locations \\
		$\mathcal{K}$ & supply types (e.g., shelter, water) \\
		$\mathcal{M}$ & transportation vehicles (e.g., truck1, truck2) \\
\end{tabular}\par}
\bigskip
{\fontsize{11}{12.5}\selectfont
\noindent\begin{tabular}{l l l}
    \multicolumn{3}{l}{\bfseries \underline{Parameters}} \\
		{\bf Symbol} & {\bf Definition} & [units] \\
		\hline
		$d_{jk}$ & demand of supply type $k$ at demand location $j$ & [units] \\
		$p_i$ & probability that pre-stage location $i$ becomes inaccessible due to the event & [-] \\
		$\tau_{ijm}$ & travel time from location $i$ to $j$ using mode $m$ & [time periods] \\
        $u_{ik}$ & pre-stage location $i$'s capacity for supply type $k$ & [units] \\
		$\alpha_k$ & weight that penalizes supply type $k$ at all pre-stage locations  & [-] \\
		$\rho$ & a sufficiently large penalty for unmet demand & [-] \\
\end{tabular}\par}
\bigskip

{\fontsize{11}{12.5}\selectfont
\noindent\begin{tabular}{l l l}
    \multicolumn{3}{l}{\bfseries \underline{Variables}} \\
		{\bf Symbol} & {\bf Definition} & [units] \\
		\hline
		$X_{ijkm}$ & quantity of supply type $k$ staged at $i$ and shipped to $j$ using mode $m$ & [units] \\
		$\bar{D}_{jk}$ & unmet demand of supply type $k$ at demand location $j$ & [units] \\
\end{tabular}\par}
\bigskip


\par\noindent{\fontsize{11}{12.5}\selectfont\bfseries\underline{Objective Function}}\par
\vspace{2pt} 

\[
\min\;
\sum_{i\in \mathcal{I}}
\sum_{j\in \mathcal{J}}
\sum_{k\in \mathcal{K}}
\sum_{m\in \mathcal{M}}
\left(
    \tau_{ijm}\,X_{ijkm}
    +
    \alpha_k\,p_i\,X_{ijkm} 
\right)
\;+\;
\sum_{j\in \mathcal{J}}
\sum_{k\in \mathcal{K}}
\rho \,\bar{D}_{jk}
\]

\par\noindent{\fontsize{11}{12.5}\selectfont\bfseries\underline{Constraints}}\par
\vspace{2pt}

\begin{subequations}\label{con:block}
\begin{align}
\sum_{j\in \mathcal{J}} \sum_{m\in \mathcal{M}} X_{ijkm}
&\le u_{ik}
&&\forall\, i\in \mathcal{I},\ k\in \mathcal{K}
&&\text{(Pre-stage capacity)}
\label{con:capacity} \\[6pt]
\sum_{i\in \mathcal{I}} \sum_{m\in \mathcal{M}} X_{ijkm}
\;+\; \bar{D}_{jk}
&\ge d_{jk}
&&\forall\, j\in \mathcal{J},\ k\in \mathcal{K}
&&\text{(Demand coverage)}
\label{con:demand} \\[6pt]
X_{ijkm} &\ge 0
&&\forall\, i\in \mathcal{I},\ j\in \mathcal{J},\ k\in \mathcal{K},\ m\in \mathcal{M}
&&\text{(Flow nonnegativity)}
\label{con:flow_nonneg} \\[4pt]
\bar{D}_{jk} &\ge 0
&&\forall\, j\in \mathcal{J},\ k\in \mathcal{K}
&&\text{(Shortage nonnegativity)}
\label{con:short_nonneg}
\end{align}
\end{subequations}

\bigskip
\clearpage

\par\noindent{\fontsize{11}{12.5}\selectfont\bfseries\underline{Constraint Description}}\par
\vspace{2pt} 

\begin{itemize}
  \item \textbf{Pre-stage capacity:}  
        Limits the total amount of each supply type shipped out of a pre-stage site to the available inventory allocated there. This is necessary when calculating the expectation metrics in post-processing.
  \item \textbf{Demand coverage:}  
        Ensures each demand location’s requirement is met either by delivered supply or, if insufficient, by recorded unmet demand.
  \item \textbf{Nonnegativity of flows:}  
        Prevents negative shipment quantities.
  \item \textbf{Nonnegativity of unmet demand:}  
        Prevents negative unmet-demand values.
\end{itemize}
\bigskip

\section{\large{Phase 1 Probabilistic Formulation, ($\mathcal{P}_{\!1}$) Post-processing: Scenario Evaluation and Expected-impact Computation}}
\bigskip

{\fontsize{11}{12.5}\selectfont
\noindent\begin{tabular}{l l}
		\multicolumn{2}{l}{\bf \underline{Sets}} \\
		{\bf Symbol} & {\bf Definition} \\
		\hline
		$\mathcal{S}$ & forecast scenarios (each with a different forecast-cone width) \\
		$\mathcal{O}$ & landfall outcomes \\
		$\mathcal{L}$ & intermediate storm locations (positions the storm may occupy before landfall) \\
		$K$ & supply types \\
		$\mathcal{T}$ & possible dispatch times \\
\end{tabular}\par}
\bigskip

\subsection{Outcome-level performance under each scenario}

\noindent
For each scenario $s \in \mathcal{S}$ and each feasible landfall outcome $o \in \mathcal{O}$, the post-processing step re-solves ($\mathcal{P}_{\!1}$) and the delivery module with:
\[
X_{ijkm}\ 
\]
where $i$ and $k$ are fixed, yielding:
\[
c_{kos} : \text{closing time for supply type } k \text{ under outcome } o \text{ in scenario } s
\]
\[
\check{d}_{kos} : \text{unmet demand for supply type } k \text{ under outcome } o \text{ in scenario } s
\]

\subsection{Conditional probabilities}

For each intermediate location $\ell \in \mathcal{L}$ and scenario $s \in \mathcal{S}$:

\[
\check{p}_{o\mid \ell,s} : \text{probability of landfall outcome } o
\text{ conditional on storm location } \ell \text{ and scenario } s
\]

\[
\hat{p}_{\ell\mid t} : \text{probability that the storm is at location } \ell
\text{ conditional on dispatch time } t
\]

\subsection{Expected performance at each intermediate location}

For each supply type $k$ and each intermediate location $\ell$:
\[
e_{k\ell}
=
\sum_{s \in \mathcal{S}} \;\sum_{o \in \mathcal{O}}
\check{p}_{o\mid \ell,s}\, c_{k o s}
\quad
\forall k \in \mathcal{K},\ \forall \ell \in \mathcal{L}
\]
\[
\bar{d}_{k\ell}
=
\sum_{s \in \mathcal{S}} \;\sum_{o \in \mathcal{O}}
\check{p}_{o\mid \ell,s}\, \check{d}_{k o s}
\quad
\forall k \in \mathcal{K},\ \forall \ell \in \mathcal{L}
\]

\subsection{Dispatch-time expected impacts}

Dispatch-time expectations are computed by aggregating over intermediate event locations for each dispatch time $t$ and supply type $k$:
\[
\hat{e}_{kt}
=
\sum_{\ell \in \mathcal{L}} \hat{p}_{\ell\mid t}\, e_{k l},
\quad
\forall k \in \mathcal{K},\ \forall t \in \mathcal{T}
\]
\[
\hat{d}_{kt}
=
\sum_{\ell \in \mathcal{L}} \hat{p}_{\ell\mid t}\, \bar{d}_{k l}
\quad
\forall k \in \mathcal{K},\ \forall t \in \mathcal{T}
\]

\bigskip

The quantities $\hat{e}_{kt}$ and $\hat{d}_{kt}$ serve as parameters in the ($\mathcal{P}_{\!2}$) optimization model.

\bigskip

\section{\large{Phase 2 Probabilistic Formulation Dispatch Timing Optimization, ($\mathcal{P}_{\!2}$)}}
\bigskip

{\fontsize{11}{12.5}\selectfont
\noindent\begin{tabular}{l l l}
		\multicolumn{3}{l}{\bf \underline{Parameters}} \\
		{\bf Symbol} & {\bf Definition} & [units] \\
		\hline
		$a_k$ & weight on expected closing time for supply type $k$ & [-] \\
		$b_k$ & weight on expected unmet demand for supply type $k$ & [-] \\
		$\hat{e}_{kt}$ & expected closing time if supply type $k$ is dispatched at time $t$ & [time] \\
		$\hat{d}_{kt}$ & expected unmet demand if supply type $k$ is dispatched at time $t$ & [\%] \\
		$\hat{\tau}_k$ & dispatch travel time from the source to pre-stage locations for supply type $k$ & [time] \\
		$\tau_k^+$ & additional dispatch-delay buffer for supply type $k$ (optionally zero) & [time] \\
\end{tabular}\par}
\bigskip

{\fontsize{11}{12.5}\selectfont
\noindent\begin{tabular}{l l l}
		\multicolumn{3}{l}{\bf \underline{Variables}} \\
		{\bf Symbol} & {\bf Definition} & [units] \\
		\hline
		$Z_{kt}$ & 1 if supply type $k$ is dispatched at time $t$, 0 otherwise & [-] \\
\end{tabular}\par}

\bigskip

\bigskip

\FloatBarrier

\par\noindent{\fontsize{11}{12.5}\selectfont\bfseries\underline{Objective Function}}\par
\vspace{2pt} 

\[
\min\;
\sum_{k\in \mathcal{K}} \sum_{t\in \mathcal{T}}
\left[
a_k\Big( \hat{e}_{kt}
+ \max\{0,\; (\hat{\tau}_k + \tau_k^+) - t\} \Big)
\;+\;
b_k\, \hat{d}_{kt}
\right] Z_{kt}
\]

\par\noindent{\fontsize{11}{12.5}\selectfont\bfseries\underline{Constraints}}\par
\vspace{2pt} 

\begin{subequations}\label{con:dispatch_block}
\begin{align}
\sum_{t\in \mathcal{T}} Z_{kt}
&= 1
&&\forall\, k\in \mathcal{K}
&&\textnormal{(Select one dispatch time)}
\label{con:select_one} \\[6pt]
Z_{kt}
&\in \{0,1\}
&&\forall\, k\in \mathcal{K},\ t\in \mathcal{T}
&&\textnormal{(Binary dispatch choice)}
\label{con:binary_dispatch}
\end{align}
\end{subequations}

\par\noindent{\fontsize{11}{12.5}\selectfont\bfseries\underline{Constraint Description}}\par
\vspace{2pt} 

\begin{itemize}
  \item \textbf{Select one dispatch time:}  
        Requires exactly one dispatch time to be chosen for each supply type.

  \item \textbf{Binary dispatch choice:}  
        Restricts the dispatch-time decision to a yes/no selection.
\end{itemize}

\bigskip

\section{\large{Phase 2 Linear Alternative ($\mathcal{P}_{\!2}$) (not used in case study):}}
\bigskip

{\fontsize{11}{12.5}\selectfont
\noindent\begin{tabular}{l l l}
		\multicolumn{3}{l}{\bf \underline{New Variable}} \\
		{\bf Symbol} & {\bf Definition} & [units] \\
		\hline
		$W_{kt}$ & auxiliary variable representing $\max\{0,\;(\hat{\tau}_k+\tau_k^+) - t\}$ & [time] \\
\end{tabular}\par}

\bigskip
\bigskip

\FloatBarrier

\par\noindent{\fontsize{11}{12.5}\selectfont\bfseries\underline{Modified Objective Function}}\par
\vspace{2pt} 

\[
\min\;
\sum_{k\in \mathcal{K}} \sum_{t\in \mathcal{T}}
\left[
a_k\Big( \hat{e}_{kt}
+ W_{kt} \Big)
\;+\;
b_k\, \hat{d}_{kt}
\right] Z_{kt}
\]


\par\noindent{\fontsize{11}{12.5}\selectfont\bfseries\underline{Modified Constraints}}\par
\vspace{2pt}

\begin{subequations}\label{con:lateness_block}
\begin{align}
W_{kt}
&\ge (\hat{\tau}_k + \tau_k^+) - t
&&\forall\, k\in \mathcal{K},\ t\in \mathcal{T}
&&\textnormal{(Lateness definition)}
\label{con:lateness_def} \\[6pt]
W_{kt}
&\ge 0
&&\forall\, k\in \mathcal{K},\ t\in \mathcal{T}
&&\textnormal{(Nonnegative lateness)}
\label{con:lateness_nonneg}
\end{align}
\end{subequations}

\par\noindent{\fontsize{11}{12.5}\selectfont\bfseries\underline{Constraint Description}}\par
\vspace{2pt}

\begin{itemize}
  \item \textbf{Lateness definition:}\normalfont\ 
        Defines the lateness amount as the shortfall between the required lead-time buffer and the chosen dispatch time.

  \item \textbf{Nonnegative lateness:}\normalfont\ 
        Prevents the lateness variable from taking negative values.
\end{itemize}
\bigskip

\section{\large{Robust Blue Formulation, ($\mathcal{R}^{+}$)}}
\bigskip

{\fontsize{11}{12.5}\selectfont
\noindent\begin{tabular}{l l}
    \multicolumn{2}{l}{\bfseries \underline{Sets}} \\
{\bf Symbol} & {\bf Definition}\\\hline
$\mathcal C,\ \mathcal V,\ \mathcal L,\ \mathcal R,\ \mathcal T,\ \mathcal G$ & commodities, individual vehicles, locations, routes, time periods, groups\\
$\bar{\mathcal R}$           & all eligible $(v,\ell,\ell',r,t)$ combinations\\
$\hat{\mathcal R}_{\ell}$ & routes with $\ell$ as the origin\\ 
$\bar{\mathcal C}$  & all eligible $(c,v,\ell,\ell',r,t)$ combinations\\
$\bar{\mathcal L}$ & all eligible $(c, \ell, t)$ combinations\\
$\mathcal E$ & events $e$, each linked to a subset of routes $(\ell,\ell',r)$

\end{tabular}\par}

\bigskip

{\fontsize{11}{12.5}\selectfont
\noindent\begin{tabular}{l l l}
    \multicolumn{3}{l}{\bfseries \underline{Parameters}} \\
		{\bf Symbol} & {\bf Definition} & [units] \\ 
		\hline
		$d_{c\ell t}$ & demand for commodity $c$ at location $\ell$ at the start of time $t$ & [units/time period] \\
		$c_e$ & impact coefficient for event $e$ & [-] \\        
        $\bar{c}_\ell$ & max vehicle departures allowed from location $\ell$ during any time period & [vehicles/time period] \\
		$s_{c\ell t}$ & external supply of commodity $c$ at location $\ell$ at the beginning of time $t$ & [units] \\
		$n_{v\ell t}$ & number of vehicles $v$ introduced at location $\ell$ and time $t$ & [vehicles] \\
		$\kappa_{v\ell\ell'r t}$ & cost of vehicle $v$ moving from $\ell$ to $\ell'$ on route $r$ at time $t$ & [cost/vehicle] \\

        $\hat{\rho}_{c\ell t}$ & penalty for unmet demand of commodity $c$ at location $\ell$ by time $t$ & [cost/unit] \\
		$\tau_{v\ell\ell'r}$ & travel time for vehicle $v$ to go from $\ell$ to $\ell'$ on route $r$ & [time periods] \\
		$\tau^{\text{turn}}_{v\ell}$ & turnaround time for vehicle $v$ at location $\ell$ & [time periods] \\
		$w_c$ & weight of commodity $c$ & [weight/unit] \\
		$\bar{w}_v$ & weight capacity of vehicle $v$ & [weight] \\
		$v_c$ & volume of commodity $c$ & [volume/unit] \\
		$\bar{v}_v$ & volume capacity of vehicle $v$ & [volume] \\

		$\theta$ & penalty coefficient for attack exposure & [cost] \\
\end{tabular}\par}
\bigskip

{\fontsize{11}{12.5}\selectfont
\noindent\begin{tabular}{l l l}
		\multicolumn{3}{l}{\bf \underline{Variables}} \\ 
		{\bf Symbol} & {\bf Definition} & [units] \\ 
		\hline
		$N_{v\ell\ell'r t}$ & 1 if vehicle $v$ departs location $\ell$ to $\ell'$ on route $r$ during time $t$, 0 otherwise & [-] \\
		$\bar{N}_{v\ell t}$ & number of vehicles $v$ remaining at location $\ell$ at the end of time $t$ & [vehicles] \\
		$X_{c v\ell\ell'r t}$ & amount of commodity $c$ traveling from $\ell$ to $\ell'$ on vehicle $v$ via route $r$  & [units] \\
        {}& that started from $\ell$ at time $t$ & {}\\
		$I_{c\ell t}$ & amount of commodity $c$ at location $\ell$ at the end of time $t$ & [units] \\
		$D_{c\ell t}$ & amount of demand for commodity $c$ satisfied at location $\ell$ during time $t$ & [units] \\
		$\bar{D}_{c\ell t}$ & unmet demand of commodity $c$ at location $\ell$ by the end of time $t$ & [units] \\
		$A$ & attack exposure variable (penalizes movement on attacked arcs) & [-] \\
\end{tabular}\par}
\bigskip

\par\noindent{\fontsize{11}{12.5}\selectfont\bfseries\underline{Objective Function}}\par
\vspace{2pt} 
\[
\min\;
\theta\,A
+\sum_{(v,\ell,\ell',r,t)\in\bar{\mathcal R}}
   \kappa_{v\ell\ell'r t}\,N_{v\ell\ell'r t}
+\sum_{(c,\ell,t)\in \bar{\mathcal L}}
   \hat\rho_{c\ell t}\,\bar D_{c\ell t}
\]

\clearpage
\par\noindent{\fontsize{11}{12.5}\selectfont\bfseries\underline{Constraints}}\par
\vspace{2pt} 

\begin{subequations}
\noindent\textbf{\bf Vehicle stock} 
\begin{align}
\bar N_{v\ell t} &=
\bar N_{v\ell,t-1}
- \sum_{\{\ell',r\} \in \bar{\mathcal R}} N_{v\ell\ell'r t}
+ \sum_{\{\ell',r\} \in \bar{\mathcal R}} 
    N_{v\ell'\ell r,\,t - \tau_{v\ell'\ell r} - \tau^{\text{turn}}_{v\ell} - 1}
+ n_{v\ell t}
&&\forall\,v \in \mathcal{V}, \ell \in \mathcal{L}, t \in \mathcal{T}
\end{align}


\noindent\textbf{\bf Commodity stock} 
\begin{align}
I_{c\ell t} &=
I_{c\ell,t-1}
- \sum_{\{v,\ell',r\} \in \bar{\mathcal C}} X_{c v \ell \ell' r t}
+ \sum_{\{v,\ell',r\} \in \bar{\mathcal C}} X_{c v \ell' \ell r,\;t - \tau_{v\ell'\ell r} - 1} \notag \\
&\quad -  \bar D_{c\ell,t-1} + \bar D_{c\ell t} - d_{c\ell t}
+ s_{c\ell t}
&& \forall (c,\ \ell,\ t) \in \bar{\mathcal L}
\end{align}

\noindent\textbf{\bf Demand inequality} 
\begin{align}
\bar D_{c\ell t} &= \bar D_{c\ell,t-1} -  D_{c\ell t} + d_{c\ell t} 
&&\forall (c, \ell,\ t) \in \bar{\mathcal L}
\end{align}

\noindent\textbf{\bf Deliver all demand} 
\[
\sum_{(c, v, \ell, \ell', r, t) \in \bar{\mathcal C}} 
    X_{c v \ell \ell' r t}
\;\ge\;
\sum_{(c, \ell', t) \in \bar{\mathcal L}} 
    d_{c\ell' t}
\quad 
\]

\noindent\textbf{\bf Weight and volume knapsacks} 
\begin{align}
\sum_{(c) \in \bar{\mathcal C}} w_c X_{c v \ell \ell' r t}
&\le \bar w_v N_{v\ell\ell' r t}
&&\forall (v,\ell,\ell',r,t) \in \bar{\mathcal R}\\
\sum_{(c) \in \bar{\mathcal C}} v_c X_{c v \ell \ell' r t}
&\le \bar v_v N_{v\ell\ell' r t}
&&\forall (v,\ell,\ell',r,t) \in \bar{\mathcal R}
\end{align}


\noindent\textbf{\bf Single trip per vehicle} 
\[
\sum_{\{\ell,\ell',r,t\} \in \bar{\mathcal R}} N_{v\ell\ell'r t} \le 1
\quad \forall v \in \mathcal{V}
\]

\noindent\textbf{\bf Attack exposure} 
\[
{A} \ge
\sum_{\substack{v,t :\\ (v,\ell,\ell',r,t) \in \bar{\mathcal R} }} 
c_e \cdot N_{v\ell\ell'r t}
\quad \forall e \in \mathcal{E}
\]

\end{subequations}

For readability, we simplify summation notation throughout model ($\mathcal{R}^{+}$) by suppressing indexed tuples and instead using $\{\cdots\}$ with the indices requiring control, e.g., $\{\ell', r\} \in \mathcal{\bar{R}}$ instead of $\ell', r: (v,\ell,\ell',r,t) \in \mathcal{\bar{R}}$.  
\bigskip
\clearpage

\par\noindent{\fontsize{11}{12.5}\selectfont\bfseries\underline{Constraint Description}}\par
\vspace{2pt} 

\begin{itemize}
  \item \textbf{Vehicle stock:}  
        Updates end-of-period vehicle availability at each location by subtracting departures during period \(t\), adding arrivals that finish travel and turnaround by the beginning of \(t\), and including any newly introduced vehicles.

  \item \textbf{Commodity stock:}  
        Balances end-of-period inventory for every commodity and location by deducting outbound shipments and demand fulfilled in period \(t\), adding inbound shipments that arrive at the beginning of \(t\), and adding external supply released at the start of \(t\).

  \item \textbf{Demand inequality:}  
        Bounds current unmet demand {\it at each location} so it can grow only by the amount of new demand issued in period \(t\), thereby precluding preemptive stockpiling.

  \item \textbf{Deliver all demand:}  
        Requires the total quantity moved across all vehicles over the planning horizon to meet or exceed the system-wide aggregate demand.

  \item \textbf{Weight knapsack:}  
        Limits the sum of commodity weights loaded on a specific vehicle trip to that vehicle’s weight capacity \(\bar w_k\).

  \item \textbf{Volume knapsack:}  
        Limits the sum of commodity volumes loaded on a specific vehicle trip to that vehicle’s volume capacity \(\bar v_k\).


  \item \textbf{Single trip per vehicle:}  
        Allows each vehicle \(k\) to initiate at most one trip in the modeled horizon. This is assuming there are enough trucks for a single dispatch. This assures that vehicles are dispatched at the beginning of the delivery operation. 

  \item \textbf{Attack exposure:}  
        Sets the exposure variable \(A\) no lower than the event-weighted sum of vehicles that traverse each edge, so its penalty can be applied in the objective. Note that there is a coefficient $c_e$ for each (origin, destination,  route) tuple; thus $Event[e]$ = $(\ell$, $\ell'$, $r$, $c_e$).
\end{itemize}
\bigskip

\section{\large{Robust Red Formulation, ($\mathcal{R}^{+}$)}}
\bigskip
{\fontsize{11}{12.5}\selectfont
\noindent\begin{tabular}{l l l}
		\multicolumn{3}{l}{\bf \underline{Parameters}} \\
		{\bf Symbol} & {\bf Definition} & [units] \\
		\hline
		$\epsilon$ & positive value for log-barrier term on individual link attack & [$\sim$arcs] \\
		$\hat{\epsilon}$ & positive value for log-barrier term on aggregate link attack & [$\sim$arcs] \\
		$n_{v\ell\ell'r t}$ & vehicle flow on each link in a given Blue solution & [units] \\
		$\beta$ & approximate total attack budget for  Red  & [cost] \\
		$b_{\ell \ell'}$ & per-link $(\ell, \ell')$ approximate attack budget & [cost/arc] \\
\end{tabular}\par}
\bigskip

{\fontsize{11}{12.5}\selectfont
\noindent\begin{tabular}{l l l}
		\multicolumn{3}{l}{\bf \underline{Variables}} \\
		{\bf Symbol} & {\bf Definition} & [units] \\
		\hline
		$W_{\ell \ell'}$ & Red's effort to disrupt Blue's flow on link $(\ell, \ell')$ & [cost/arc] \\
		$Z_{\ell \ell'}$ & per-unit cost change derived from Red's disruption on link $(\ell, \ell')$ & [cost/unit] \\
	\end{tabular}\par}
\bigskip

\clearpage
\par\noindent{\fontsize{11}{12.5}\selectfont\bfseries\underline{Objective Function}}\par
\vspace{2pt} 

\begin{equation}
(\mathcal{R}) \;\;\; \max \quad
\sum_{\substack{(\ell, \ell') \in \mathcal{\bar{R}} \\ N_{v\ell\ell'r t} > 0}}
n_{v \ell \ell'r t}Z_{\ell \ell'}
\;+\;
\epsilon \cdot
\sum_{\substack{(\ell, \ell') \in \mathcal{\bar{R}} \\ N_{v\ell\ell'r t} > 0}}
\log\left(b_{\ell \ell'} - W_{\ell \ell'}\right)
\;+\;
\hat{\epsilon} \cdot
\log\left(\beta -
\sum_{\substack{(\ell, \ell') \in \mathcal{\bar{R}} \\ N_{v\ell\ell'r t} > 0}}
W_{\ell \ell'}\right)
\label{redobjfn}
\end{equation}
\bigskip

\par\noindent{\fontsize{11}{12.5}\selectfont\bfseries\underline{Constraints}}\par
\vspace{2pt} 

\begin{subequations}\label{con:red_block}
\begin{align}
Z_{\ell \ell'}
&= e^{W_{\ell \ell'}} - 1
&&\forall\, (\ell, \ell') \in \mathcal{\bar{R}}
&&\textnormal{(Effort-to-impact mapping)}
\label{con:red_mapping} \\[6pt]
0 \;\le\; W_{\ell \ell'}
&\le 1
&&\forall\, (\ell, \ell') \in \mathcal{\bar{R}}
&&\textnormal{(Effort domain)}
\label{con:red_bounds}
\end{align}
\end{subequations}
\bigskip

\par\noindent{\fontsize{11}{12.5}\selectfont\bfseries\underline{Constraint Description}}\par
\vspace{2pt} 

\begin{itemize}
  \item \textbf{Effort-to-impact mapping:}  
        Converts Red's disruption effort on a link into the cost increase applied to that link in the next Blue solve.

  \item \textbf{Effort domain:}  
        Limits Red's disruption effort on any single link to between 0 and 1, inclusive.
\end{itemize}

\noindent\textbf{Note:} The per-link cost increase produced by Red is carried forward into the next iteration as the link impact coefficient used in the Blue formulation.

\end{APPENDIX}
%
%

\ACKNOWLEDGMENT{We would like to express our sincere gratitude to Caleb Fluker (Systems Engineering, West Point) and Dr. David Ciemnoczolowski (Defcon AI) for their technical insights and to  Dr. Jeff Grobman (Defcon AI) for providing us with this research opportunity.  }







\bibliographystyle{informs2014} 
\bibliography{sample.bib}

\end{document}